\documentclass[aap,seceqn,citesort,MSNbibl,dvips]{arximspdf}
\usepackage{graphicx}

\doi{10.1214/10-AAP679}
\volume{20}
\issue{6}
\pubyear{2010}
\firstpage{2118}
\lastpage{2161}

\makeatletter

\newtheorem{theorem}{Theorem}[section]
\newtheorem{lemma}[theorem]{Lemma}
\newtheorem{cor}[theorem]{Corollary}
\newtheorem{prop}[theorem]{Proposition}

\makeatother

\begin{document}
\begin{frontmatter}

\title{Asymptotic behavior of the
finite-size magnetization as a function
of the speed of approach to criticality}
\runtitle{Asymptotic behavior of the
finite-size magnetization}

\begin{aug}
\author[A]{\fnms{Richard S.} \snm{Ellis}\corref{}\thanksref{t1}\ead[label=e1]{rsellis@math.umass.edu}},
\author[B]{\fnms{Jonathan} \snm{Machta}\ead[label=e2]{machta@physics.umass.edu}} and
\author[C]{\fnms{Peter Tak-Hun} \snm{Otto}\ead[label=e3]{potto@willamette.edu}}
\runauthor{R. S. Ellis, J. Machta and P. T.-H. Otto}
\affiliation{University of Massachusetts, University of Massachusetts
and~Willamette~University}
\address[A]{R. S. Ellis\\
Department of Mathematics and Statistics\\
University of Massachusetts\\
Amherst, Massachusetts 01003\\
USA\\
\printead{e1}} 
\address[B]{J. Machta\\
Department of Physics\\
University of Massachusetts\\
Amherst, Massachusetts 01003\\
USA\\
\printead{e2}}
\address[C]{P. T.-H. Otto\\
Department of Mathematics\\
Willamette University\\
Salem, Oregon 97301\\
USA\\
\printead{e3}}
\end{aug}

\thankstext{t1}{Supported by NSF Grant DMS-06-04071.}

\received{\smonth{6} \syear{2009}}
\revised{\smonth{1} \syear{2010}}

%
\begin{abstract}
The main focus of this paper is to determine whether the thermodynamic
magnetization is
a physically relevant estimator of the finite-size magnetization. This
is done by comparing the asymptotic
behaviors of these two quantities along parameter
sequences converging to either a second-order point or the tricritical
point in the mean-field
Blume--Capel model. We show that the thermodynamic magnetization and
the finite-size magnetization
are asymptotic when the parameter $\alpha$ governing the speed at which
the sequence approaches criticality is below a certain threshold
$\alpha_0$. However, when $\alpha$ exceeds
$\alpha_0$, the thermodynamic magnetization converges to 0 much faster
than the finite-size magnetization.
The asymptotic behavior of the finite-size magnetization is proved via
a moderate deviation principle when
$0 < \alpha< \alpha_0$ and via a weak-convergence limit when $\alpha
> \alpha_0$.
To the best of our knowledge, our results are the first rigorous
confirmation of the
statistical mechanical theory of finite-size scaling for a mean-field model.
\end{abstract}

%

\end{frontmatter}

\section{Introduction}
\label{section:intro}

For the mean-field Blume--Capel model, as for other mean-field spin systems,
the magnetization in the thermodynamic limit is well understood within
the theory of large deviations.
In this framework the thermodynamic magnetization arises as the unique,
positive, global minimum point
of the rate function in a large deviation principle.
The question answered in this paper is whether, in a neighborhood of
criticality,
the thermodynamic magnetization is a physically relevant estimator of the
finite-size magnetization, which is the expected value of the spin per
site. A similar question is answered by
the heuristic, statistical mechanical theory of finite-size scaling.
This paper is both motivated by the theory of finite-size scaling
and puts that theory on a firm foundation in the context of mean-field
spin systems. It is hoped that our results
suggest how this question can be addressed in the context of much more
complicated, short-range spin systems.

Our approach is to evaluate the asymptotic behaviors of the
thermodynamic magnetization and the physically relevant,
finite-size magnetization along parameter
sequences converging to either a second-order point or the tricritical
point in the mean-field Blume--Capel model.
The thermodynamic magnetization is then considered to be a physically
relevant estimator of the
finite-size magnetization when these two quantities have the same
asymptotic behavior.
Our main finding is that the value of the parameter $\alpha$ governing
the speed at which the sequence approaches criticality
determines whether or not the asymptotic behaviors of these two
quantities are the same.
Specifically, we show in Theorem \ref{thm:mainsmall} that the
thermodynamic magnetization and the finite-size magnetization
are asymptotic when $\alpha$ is below a certain threshold $\alpha_0$
and that therefore
the thermodynamic magnetization is a physically relevant estimator when
$0 < \alpha< \alpha_0$. However, when $\alpha$ exceeds
$\alpha_0$, then according to Theorem \ref{thm:mainlarge},
the thermodynamic magnetization converges to 0 much faster than the
finite-size magnetization, and therefore
the thermodynamic magnetization is not a physically relevant estimator
when $\alpha> \alpha_0$. An advantage of using
the thermodynamic magnetization as an estimator of the finite-state
magnetization when $0 < \alpha< \alpha_0$
is that the asymptotic behavior of the former quantity is much easier
to derive than the
asymptotic behavior of the latter quantity [see the discussion at the end of
the paragraph after (\ref{eqn:lim5})].

The investigation is carried out for a mean-field version of an important
lattice spin model due to Blume and Capel, to which we refer as the
B--C model \cite{Blu,Cap1,Cap2,Cap3}.
This mean-field model is equivalent to the B--C model on
the complete graph on $N$ vertices. It is one of the simplest models that
exhibits the following intricate phase-transition structure: a curve of
second-order
points, a curve of first-order points and
a tricritical point, which separates the two curves.
A generalization of the B--C model is studied in \cite{BluEmeGri}.


The mean-field B--C model is defined by a canonical ensemble that we
denote by $P_{N,\beta,K}$;
$N$ equals the number of spins, $\beta$ is the inverse temperature
and $K$ is the interaction strength. $P_{N,\beta,K}$ is defined
in (\ref{eqn:pnbetak}) in terms of the Hamiltonian
\[
H_{N, K}(\omega)=\sum_{j=1}^N \omega_j^2 - \frac{K}{N} \Biggl( \sum_{j=1}^N
\omega_j \Biggr)^2,
\]
in which $\omega_j$ represents the spin at site $j \in\{1,2,\ldots,
N\}$ and
takes values in $\Lambda=\{1,0,-1\}$.
The configuration space for the model is the set $\Lambda^N$
containing all
sequences $\omega= (\omega_1,\omega_2, \ldots, \omega_N)$ with
each $\omega_j \in\Lambda$.
Expectation with respect to $P_{N,\beta,K}$ is denoted by $E_{N,\beta
,K}$. The finite-size magnetization
is defined by $E_{N,\beta,K}\{|S_N/N|\}$, where $S_N$ equals the total
spin $\sum_{j=1}^N \omega_j$.

Before introducing the results in this paper,
we summarize the phase-transition structure of the model.
For $\beta> 0$ and $K >0$ we denote by $\mathcal{M}_{\beta,K}$ the
set of
equilibrium values
of the magnetization.
$\mathcal{M}_{\beta,K}$ coincides with the
set of global minimum points of the free-energy functional $G_{\beta
,K}$, which
is defined in (\ref{eqn:introgbetak}). It is known from
heuristic arguments
and is proved in \cite{EllOttTou} that there exists a critical
inverse temperature $\beta_c = \log4$ and that for $0 < \beta\leq
\beta_c$
there exists a quantity $K(\beta)$ and for $\beta> \beta_c$ there exists
a quantity $K_1(\beta)$ having the following properties. The positive
quantity $m(\beta,K)$
appearing in the following list is the thermodynamic magnetization.
\begin{enumerate}
\item Fix $0 < \beta\leq\beta_c$. Then for $0 < K \leq K(\beta)$,
$\mathcal{M}_{\beta,K}$ consists of the unique
pure phase 0, and for $K > K(\beta)$, $\mathcal{M}_{\beta,K}$
consists of two
nonzero values $\pm m(\beta,K)$.
\item For $0 < \beta\leq\beta_c$, $\mathcal{M}_{\beta,K}$
undergoes a continuous bifurcation at $K = K(\beta)$, changing continuously
from $\{0\}$ for $K \leq K(\beta)$ to $\{\pm m(\beta,K)\}$ for $K >
K(\beta)$.
This continuous bifurcation corresponds to a second-order phase transition.
\item Fix $\beta> \beta_c$. Then for $0 < K < K_1(\beta)$, $\mathcal
{M}_{\beta,K}$
consists of the unique
pure phase 0, for $K = K_1(\beta)$, $\mathcal{M}_{\beta,K}$ consists
of 0
and two nonzero values $\pm m(\beta,K_1(\beta))$
and for $K > K_1(\beta)$, $\mathcal{M}_{\beta,K}$ consists of two
nonzero values $\pm m(\beta,K)$.
\item For $\beta> \beta_c$, $\mathcal{M}_{\beta,K}$
undergoes a discontinuous bifurcation at $K = K_1(\beta)$,
changing discontinuously from $\{0\}$ for $K < K(\beta)$ to
$\{0, \pm m(\beta,K)\}$ for $K = K_1(\beta)$ to
$\{\pm m(\beta,K)\}$ for $K > K_1(\beta)$. This discontinuous bifurcation
corresponds to a first-order phase transition.
\end{enumerate}

\begin{figure}

\includegraphics{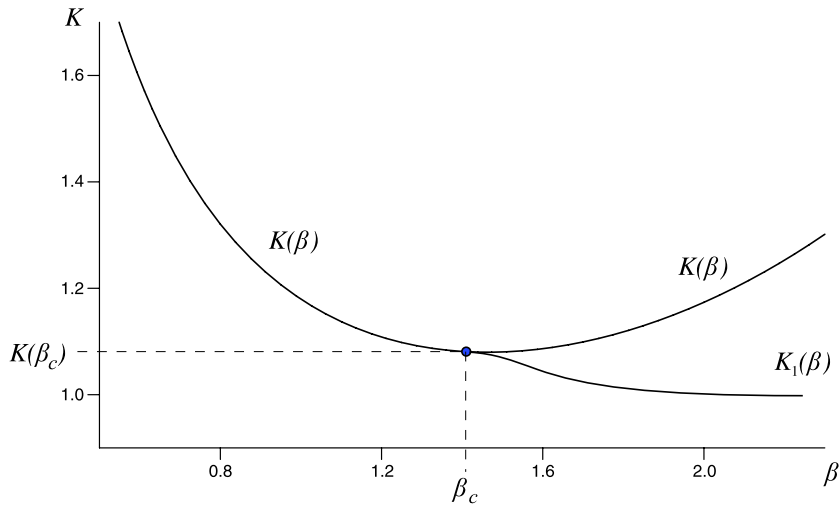}

\caption{The sets that describe the phase-transition structure of the
mean-field B--C model:
the second-order curve $\{(\beta,K(\beta)), 0 < \beta< \beta_c\}$, the
first-order curve $\{(\beta,K_1(\beta)), \beta> \beta_c\}$ and the
tricritical point $(\beta_c,K(\beta_c))$. The phase-coexistence
region consists of all $(\beta,K)$ above the second-order curve, above
the tricritical point, on the first-order
curve and above the first-order curve. The extension of the
second-order curve to $\beta> \beta_c$
is called the spinodal curve.}
\label{figure1}
\end{figure}

Because of items 2 and 4, we refer to the curve $\{(\beta,K(\beta)),
0 < \beta
< \beta_c\}$ as the second-order curve and to the curve
$\{(\beta,K_1(\beta)), \beta> \beta_c\}$ as the first-order curve.
Points on the second-order curve are called second-order points, and points
on the first-order curve first-order points.
The point $(\beta_c,K(\beta_c)) = (\log4, 3/2 \log4)$ separates the
second-order curve
from the first-order curve and is called the tricritical point. The two-phase
region consists of all points in the positive $\beta$-$K$ quadrant for which
$\mathcal{M}_{\beta,K}$ consists of two values. Thus this region
consists of all
$(\beta,K)$ above the second-order curve, above the tricritical point
and above the first-order curve; that is, all
$(\beta,K)$ satisfying $0 < \beta\leq\beta_c$ and $K > K(\beta)$ and
satisfying $\beta> \beta_c$ and $K > K_1(\beta)$. The sets that describe
the phase-transition structure
of the model are shown in Figure~\ref{figure1}.

For fixed $(\beta,K)$ lying in the two-phase region the finite-size
magnetization $E_{N,\beta,K}\{|S_N/n|\}$ converges to
the thermodynamic magnetization $m(\beta,K)$ as \mbox{$N \rightarrow\infty
$}. In
order to see this, we use the large deviation principle
(LDP) for $S_N/N$ with respect to $P_{N,\beta,K}$ in
\cite{EllOttTou}, Theorem 3.3,
and the fact that the set of global minimum points of the rate function
in that LDP
coincides with the set $\mathcal{M}_{\beta,K}$~\cite{EllOttTou}, Proposition 3.4,
the structure of which has just been described.
Since for $(\beta,K)$ lying in the two-phase region $\mathcal
{M}_{\beta,K}= \{\pm
m(\beta,K)\}$, the LDP implies
that the $P_{N,\beta,K}$-distributions of $S_N/N$ put an exponentially
small mass on the complement of any
open set containing $\pm m(\beta,K)$. Symmetry then yields the
weak-convergence limit
%
\begin{equation}
\label{eqn:lim1}
P_{N,\beta,K}\{S_N/N \in dx\} \Longrightarrow\bigl(\tfrac{1}{2}\delta
_{m(\beta,K)} + \tfrac{1}{2}\delta_{-m(\beta,K)} \bigr)(dx).
\end{equation}
This implies the desired result
%
\begin{equation}
\label{eqn:lim2}
\lim_{N \rightarrow\infty} E_{N,\beta,K}\{|S_N/N|\} = m(\beta,K).
\end{equation}
%

The limit in the last display is closely related to the main focus of
this paper. It shows that because the thermodynamic
magnetization is the limit, as the number of spins goes to $\infty$,
of the
finite-size magnetization, the thermodynamic magnetization $m(\beta
,K)$ is a physical relevant estimator
of the finite-size magnetization, at least when evaluated at fixed
$(\beta,K)$ in the two-phase region.

The main focus of this paper is to determine whether the thermodynamic
magnetization is a physically relevant estimator
of the finite-size magnetization in a more general sense, namely, when
evaluated along a class of sequences $(\beta_n,K_n)$ that
converge to a second-order point $(\beta,K(\beta))$ or the
tricritical point $(\beta_c,K(\beta_c))$.
The criterion for determining whether $m(\beta_n,K_n)$ is a physically
relevant estimator
is that as $n \rightarrow\infty$, $m(\beta_n,K_n)$ is asymptotic to the
finite-size magnetization $E_{n,\beta_n,K_n}\{|S_n/n|\}$, both of
which converge
to 0.
In this formulation we let $N = n$ in the finite-size magnetization;
that is, we let the number of spins $N$ coincide with
the index $n$ parametrizing the sequence $(\beta_n,K_n)$.
As summarized in Theorems \ref{thm:mainsmall}
and \ref{thm:mainlarge}, our main finding is that $m(\beta_n,K_n)$ is a
physically relevant estimator if
the parameter $\alpha$ governing the speed at which
$(\beta_n,K_n)$ approaches criticality is below a certain threshold
$\alpha_0$; however, this is not true if $\alpha > \alpha_0$.
For the sequences under consideration the parameter $\alpha$ determines
the limits
\[
b = \lim_{n \rightarrow\infty} n^\alpha(\beta_n- \beta)
\quad\mbox{and}\quad
k = \lim_{n \rightarrow\infty} n^\alpha\bigl(K_n- K(\beta)\bigr),
\]
which are assumed to exist and not to be both 0. The value of $\alpha_0$
depends on the type of the phase transition---first-order,
second-order or tricritical---that influences the sequence, an issue
addressed in Section 5 of \cite{EllMacOtt1}.

\begin{figure}

\includegraphics{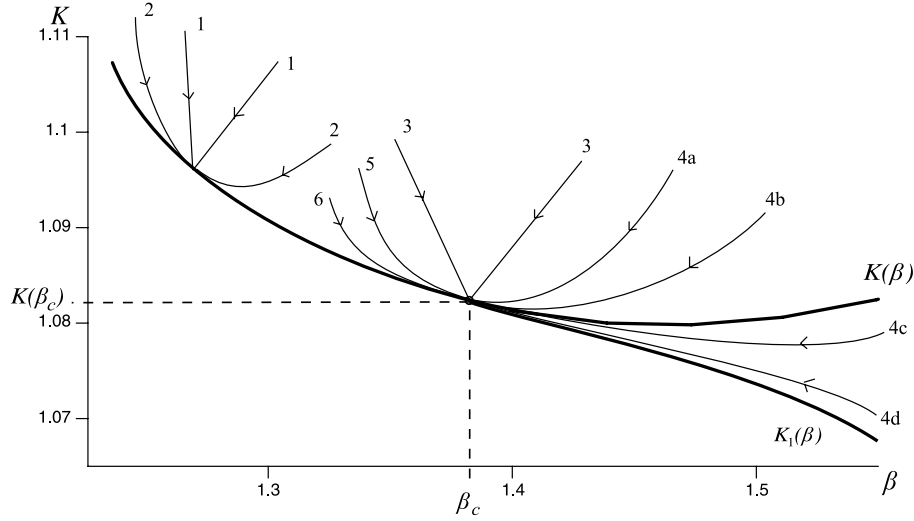}

\caption{Possible paths for the six sequences converging to a
second-order point and to the tricritical point.
The asymptotic results for the sequences converging on the paths
labeled 1, 2, 3, \textup{4a--4d}, 5 and~6 are discussed in the respective
Theorems \protect\ref{thm:resultsone}--\protect\ref{thm:resultssix}.
The sequences on the paths labeled \textup{4a--4d} are defined in
(\protect\ref{eqn:bnfour}) and in the second paragraph after that equation.}
\label{figure2}
\end{figure}

We illustrate the results contained in these two theorems by applying
them to six types of sequences.
In the case of second-order points
two such sequences are considered in Theorems \ref{thm:resultsone} and
\ref{thm:resultstwo}, and in the case of the tricritical
point four such sequences are considered in Theorems \ref
{thm:resultsthree}--\ref{thm:resultssix}.
Possible paths followed by these sequences are shown in Figure \ref{figure2}. We
believe that modulo
uninteresting scale changes, irrelevant higher order terms and other
inconsequential modifications,
these are all the sequences of the form $\beta_n = \beta+ b/n^{\alpha}$
and $K_n$ equal to $K(\beta)$ plus a polynomial in $1/n^\alpha$,
where $(\beta,K(\beta))$ is
either a second-order point or the tricritical point and $m(\beta_n,K_n)
\sim c/n^\delta$ for some $c > 0$ and $\delta> 0$.

We next summarize our main results on the asymptotic behaviors of the
thermodynamic magnetization
and the finite-size magnetization, first for small values of $\alpha$
and then for large values of $\alpha$.
The relevant information is given, respectively, in Theorems \ref
{thm:exactasymptotics},
\ref{thm:mainsmall} and \ref{thm:mainlarge}. These theorems are valid
for suitable
positive sequences $(\beta_n,K_n)$ parametrized by $\alpha> 0$, lying in
the two-phase region for all sufficiently
large $n$, and converging either to a second-order point or to the
tricritical point. The hypotheses of these three theorems overlap but
do not coincide.
The hypotheses of Theorem \ref{thm:exactasymptotics} are satisfied by
all six sequences considered
in Section~\ref{section:sequences} while the hypotheses of each of the
Theorems \ref{thm:mainsmall} and
\ref{thm:mainlarge} are satisfied by all six sequences with one exception.
For each of the six sequences the quantities $\theta$ and $\alpha_0$
appearing in these asymptotic results are specified in Table \ref{table1}.

%
%
\begin{table}
\caption{The equations where each of the six sequences is defined,
the theorems where the asymptotic results in (\protect\ref{eqn:lim3}),
(\protect\ref{eqn:lim4}) and (\protect\ref{eqn:lim5})
are stated for each sequence, and the values of $\alpha_0$, $\theta$ and
$\kappa=\frac{1}{2}(1-\alpha/\alpha_0) + \theta\alpha$
[see (\protect\ref{eqn:lim6})] for each sequence}
\label{table1}
%
\begin{tabular*}{\tablewidth}{@{\extracolsep{\fill}}lccccc@{}}
\hline
\textbf{Seq.} & \textbf{Defn.} & \textbf{Thm.} & $\bolds{\alpha_0}$
& $\bolds{\theta}$ & $\bolds{\kappa}$ \\
\hline
1 & (\ref{eqn:bnone}) & Theorem \ref{thm:resultsone} &
$\frac{1}{2}$ & $\frac{1}{2}$ & $\frac{1}{2}(1-\alpha)$ \\[4pt]
2 & (\ref{eqn:bntwo}) & Theorem \ref{thm:resultstwo} &
$\frac{1}{2p}$ & $\frac{p}{2}$ & $\frac{1}{2}(1-p\alpha)$ \\[4pt]
3 & (\ref{eqn:bnthree}) & Theorem \ref{thm:resultsthree} &
$\frac{2}{3}$ & $\frac{1}{4}$ & $\frac{1}{2}(1-\alpha)$ \\[4pt]
4 & (\ref{eqn:bnfour}) & Theorem \ref{thm:resultsfour} &
$\frac{1}{3}$ & $\frac{1}{2}$ & $\frac{1}{2}(1 - 2\alpha)$ \\[4pt]
5 & (\ref{eqn:bnfive}) & Theorem \ref{thm:resultsfive} &
$\frac{1}{3}$ & $\frac{1}{2}$ & $\frac{1}{2}(1- 2\alpha)$ \\[4pt]
6 & (\ref{eqn:bnsix}) & Theorem \ref{thm:resultssix} &
$\frac{1}{2p-1}$ & $\frac{1}{2}(p-1)$ & $\frac{1}{2}(1-p\alpha)$ \\
\hline
\end{tabular*}
\end{table}

The difference in the asymptotic behaviors of the thermodynamic
magnetization and the finite-size magnetization
for $\alpha> \alpha_0$ is described in item 3. As we discuss in
Section \ref{section:fss}, the difference is explained
by the statistical mechanical theory of finite-size scaling.

\begin{enumerate}
\item According to Theorem \ref{thm:exactasymptotics},
there exists positive quantities $\bar{x}$ and $\theta$ such that for
all $\alpha> 0$
%
\begin{equation}
\label{eqn:lim3}
m(\beta_n,K_n)\sim\bar{x}/n^{\theta\alpha}.
\end{equation}
\item ($0 < \alpha< \alpha_0$). According to
Theorem \ref{thm:mainsmall}, there exists a threshold
value $\alpha_0 > 0$ such that for all $0 < \alpha< \alpha_0$
%
\begin{equation}
\label{eqn:lim4}
E_{n,\beta_n, K_n}\{|S_n/n|\} \sim\bar{x}/n^{\theta\alpha}
\quad\mbox{and}\quad
E_{n,\beta_n, K_n}\{|S_n/n|\} \sim m(\beta_n, K_n).
\end{equation}
Because $m(\beta_n,K_n)$ is asymptotic to the finite-size magnetization,
$m(\beta_n,K_n)$ is a physically relevant estimator
of the finite-size magnetization.
In this case $(\beta_n,K_n)$ converges to criticality slowly, and we are
in the two-phase region, where the system is effectively infinite.
Formally the first index $n$ parametrizing the finite-size
magnetization can be sent to $\infty$ before the index
$n$ parametrizing the sequence $(\beta_n,K_n)$ is sent to $\infty$, and
so we have
\[
E_{n,\beta_n, K_n}\{|S_n/n|\} \approx\lim_{N \rightarrow\infty}
E_{N,\beta_n, K_n}\{|S_N/N|\} = m(\beta_n,K_n).
\]
\item ($\alpha> \alpha_0$). According to
Theorem \ref{thm:mainlarge}, there
exists a positive quantity $\bar{y}$ such that for all $\alpha>
\alpha_0$
%
\begin{eqnarray}
\label{eqn:lim5}
E_{n,\beta_n, K_n}\{|S_n/n|\} &\sim& \bar{y}/n^{\theta\alpha_0}\quad
\mbox{and}\nonumber\\[-8pt]\\[-8pt]
E_{n,\beta_n, K_n}\{|S_n/n|\} &\gg& m(\beta_n, K_n) \sim\bar{x}
/n^{\theta\alpha}.\nonumber
\end{eqnarray}
Because $m(\beta_n,K_n)$ converges to 0 much faster than the finite-size
magnetization, $m(\beta_n,K_n)$ is not a physically relevant
estimator of the finite-size magnetization.
In this case $(\beta_n,K_n)$ converges to criticality quickly, and we are
in the critical regime, where finite-size scaling effects
are important.
\end{enumerate}

The asymptotic behavior of the thermodynamic magnetization $m(\beta_n,K_n)
\rightarrow0$
stated in (\ref{eqn:lim3}) holds for all $\alpha> 0$.
It is derived in Theorem 3.2 in \cite{EllMacOtt1} and is summarized in
Theorem \ref{thm:exactasymptotics} in the present paper.
In (\ref{eqn:lim4}) we state the asymptotic behavior of the
finite-size magnetization $E_{n,\beta_n,K_n}\{|S_n/n|\} \rightarrow0$
for $0 < \alpha< \alpha_0$. This result is proved in part (a) of Theorem
\ref{thm:mainsmall} as a consequence of the moderate deviation
principle (MDP) for the spin in Theorem \ref{thm:mdp},
the weak-convergence limit in Corollary \ref{cor:weaklimit}, and the
uniform integrability estimate in Lemma \ref{lem:uniformint}.
The asymptotic behavior of $E_{n,\beta_n,K_n}\{|S_n/n|\} \rightarrow
0$ stated in
(\ref{eqn:lim5}) for $\alpha> \alpha_0$
is proved in part (a) of Theorem \ref{thm:mainlarge}
as a consequence of the weak-convergence limit for the spin in Theorem
\ref{thm:scalinglargealpha} and the
uniform-integrability-type estimate in Proposition \ref{prop:weakerunifint}.
In part (a) of Theorem \ref{thm:mainmiddle} we state the asymptotic
behavior of
$E_{n,\beta_n,K_n}\{|S_n/n|\} \rightarrow0$ for $\alpha= \alpha_0$.
That result is a consequence of a weak-convergence limit analogous to
the limit in Theorem \ref{thm:scalinglargealpha}
and the uniform-integrability-type estimate in Proposition \ref
{prop:weakerunifint}.
With changes in notation only, Theorem \ref{thm:exactasymptotics} and
Theorems \ref{thm:mainsmall}--\ref{thm:mainmiddle}
also apply to other mean-field models including the Curie--Weiss model
\cite{Ellis} and the Curie--Weiss--Potts model \cite{EllWan}.
The proof of the asymptotic behavior of the thermodynamic magnetization
in \cite{EllMacOtt1}, Theorem 3.2, is purely analytic
and is much more straightforward than the probabilistic proofs of the
asymptotic behaviors of the finite-size magnetization in
Theorems \ref{thm:mainsmall}--\ref{thm:mainmiddle}.

\begin{figure}[b]

\includegraphics{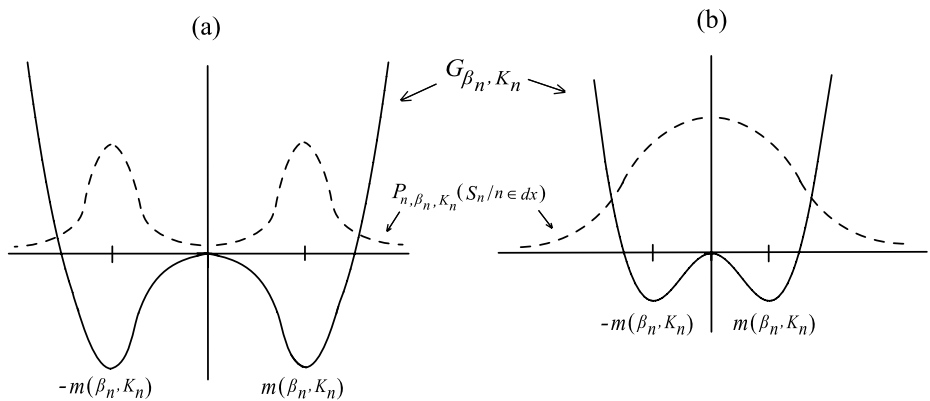}

\caption{$G_{\beta_n,K_n}$ and $P_{n,\beta_n,K_n}\{S_n/n \in dx\}$
for \textup{(a)}
$0 < \alpha< \alpha_0$, \textup{(b)} $\alpha> \alpha_0$.
Graph \textup{(b)} is not shown to scale. In fact, for $\alpha> \alpha_0$ the
global minimum points $\pm m(\beta_n,K_n)$ of $G_{\beta_n,K_n}$ are
much closer
to the origin and are much shallower than shown in graph \textup{(b)}.}
\label{figure3}
\end{figure}

Figure \ref{figure3} gives a pictorial representation of the phenomena that are
summarized in
(\ref{eqn:lim4}) for $0 < \alpha< \alpha_0$ and in
(\ref{eqn:lim5}) for $\alpha> \alpha_0$.
As we discuss in Section~\ref{section:phasetr}, for the sequences
$(\beta_n,K_n)$ under consideration
the thermodynamic magnetization $m(\beta_n,K_n)$ can be characterized
as the
unique, positive, global minimum
point in an LDP or, equivalently,
as the unique, positive, global minimum point of the dual, free-energy
functional $G_{\beta_n,K_n}$ defined in (\ref{eqn:introgbetak}).
According to graph (a) in Figure~\ref{figure3}, for $0 < \alpha< \alpha_0$,
$G_{\beta_n,K_n}$ has two deep, global minimum points at $\pm m(\beta_n,K_n)$.
Graph (b) in Figure \ref{figure3}, which is not shown to scale, exhibits the
contrasting situation for $\alpha> \alpha_0$.
In this case the global minimum points of $G_{\beta_n,K_n}$ at $\pm
m(\beta_n,K_n)$
are shallow
and close to the origin.
In the two graphs we also show the form of the distribution $P_{n,\beta
_n,K_n}\{
S_n/n \in dx\}$.
For $0 < \alpha< \alpha_0$ this probability distribution is sharply
peaked at $\pm m(\beta_n,K_n)$ as $n \rightarrow\infty$. In contrast,
for $\alpha> \alpha_0$ the probability distribution is peaked at 0
and its standard deviation is much larger than $m(\beta_n,K_n)$.

In a work in progress
we refine the asymptotic result in (\ref{eqn:lim4}), which
states that for $0 < \alpha< \alpha_0$, $m(\beta_n,K_n)$ is asymptotic
to $E_{n,\beta_n,K_n}\{|S_n/n|\}$ as $n \rightarrow\infty$. Define
$\kappa=\frac{1}{2} (1 -
\alpha/\alpha_0) + \theta\alpha$, which exceeds $\theta\alpha$
since $1 - \alpha/\alpha_0 > 0$. We conjecture that for a class of
suitable sequences $(\beta_n,K_n)$ that includes
the first five sequences considered in Section \ref{section:sequences},
there exists a positive quantity $\bar{v}$ such that for all \mbox{$0 <
\alpha< \alpha_0$}
%
\begin{equation}
\label{eqn:lim6}
E_{n,\beta_n, K_n}\bigl\{\bigl| |S_n/n| - m(\beta_n, K_n)\bigr|\bigr\} \sim\bar
{v}/n^\kappa.
\end{equation}
This refined asymptotic result would extend part (b) of Theorem \ref
{thm:mainsmall}. It is a consequence of the conjecture that
when $S_n/n$ is conditioned to lie in a suitable neighborhood of
$m(\beta_n,K_n)$, the
$P_{n,\beta_n,K_n}$-distributions of $n^\kappa(S_n/n - m(\beta
_n,K_n))$ converge in
distribution to a Gaussian.


For easy reference we list in Table \ref{table1} information about the six
sequences considered in
Section \ref{section:sequences}. The first two columns list,
respectively, the equation in which
each sequence is defined and the theorem in which the asymptotic
results in equations (\ref{eqn:lim3}),
(\ref{eqn:lim4}) and (\ref{eqn:lim5}) are stated for each sequence.
In these theorems the quantities $\bar{x}$ and $\bar{y}$ appearing in the
three asymptotic results are defined.
The three asymptotic results involve the quantities
$\alpha_0$, $\theta$, $\theta\alpha$ and $\theta\alpha_0$, the
values of the first two of which are listed in the next
two columns of the table.
In the last column of the table we list the values
of $\kappa=\frac{1}{2} (1 - \alpha/\alpha_0) + \theta\alpha$.
Through the factor $n^{-\kappa}$, $\kappa$
governs the conjectured asymptotics of $E_{n,\beta_n,K_n}\{| |S_n/n| -
m(\beta_n,K_n)|\}
$ stated in~(\ref{eqn:lim6}).

The contents of this paper are as follows. In Section
\ref{section:phasetr} we summarize the phase-transition
structure of the mean-field B--C model. Theorem \ref{thm:exactasymptotics}
in Section \ref{section:mbnkn} gives the asymptotic behavior of the
thermodynamic magnetization
$m(\beta_n,K_n)\rightarrow0$ for suitable sequences $(\beta_n,K_n)$
converging either
to a second-order point or to the tricritical point.
The heart of the paper is Section \ref{section:rates}. In this section
Theorems~\ref{thm:mainsmall},
\ref{thm:mainlarge} and \ref{thm:mainmiddle} give the asymptotic
behavior of the finite-size magnetization $E_{n,\beta_n,K_n}\{
|S_n/n|\} \rightarrow0$
for three respective ranges of $\alpha\dvtx 0 < \alpha< \alpha_0$,
$\alpha> \alpha_0$ and $\alpha= \alpha_0$. The quantity $\alpha_0$
is a threshold value that
depends on the type of the phase transition---first-order,
second-order or tricritical---that
influences the associated sequence $(\beta_n,K_n)$. These theorems also
compare the asymptotic
behaviors of the thermodynamic magnetization and the finite-size magnetization,
showing that they are the same for $0 < \alpha < \alpha_0$ but not the
same for $\alpha>\alpha_0$.
In Section \ref{section:sequences} the three theorems in the preceding
section are applied to six
specific sequences $(\beta_n,K_n)$, the first two of which converge to a
second-order point and the last
four of which converge to the tricritical point. Section
\ref{section:fss} gives\vadjust{\goodbreak} an overview of the statistical
mechanical theory of finite-size scaling, which gives insight into the
physical phenomena underlying our mathematical results.
Part (a) of Theorem \ref{thm:mainsmall} is derived in Section~\ref
{section:smallalpha}
from the MDP for the spin in Theorem \ref{thm:mdp},
the weak-convergence limit for the spin in Corollary \ref{cor:weaklimit},
and the uniform integrability estimate in Lemma \ref{lem:uniformint}.
Finally, part (a) of Theorem \ref{thm:mainlarge} is derived in
Section \ref{section:largealpha} from the weak-convergence limit for
the spin
in Theorem \ref{thm:scalinglargealpha} and the uniform-integrability-type
estimate in Proposition \ref{prop:weakerunifint}.


\section{Phase-transition structure of the mean-field B--C model}
\label{section:phasetr}

After defining the mean-field B--C model, we introduce a function
$G_{\beta,K}
$, called the free-energy functional.
The global minimum points of this function define the equilibrium
values of the magnetization.
The phase-transition structure of the model is summarized in Theorems
\ref{thm:secondorder} and \ref{thm:firstorder}. The first theorem
shows that the model exhibits
a second-order phase transition for $\beta\in(0,\beta_c]$, where
$\beta_c=
\log4$ is the critical
inverse temperature of the model. The second theorem shows that the
model exhibits a first-order
phase transition for $\beta> \beta_c$.

For $N \in\mathbb N$ the mean-field B--C model is a lattice-spin model
defined on the complete graph on $N$ vertices $1,2, \ldots, N$.
The spin at site $j \in\{1,2,\ldots, N\}$ is denoted by $\omega_j$,
a quantity taking values in $\Lambda=\{1,0,-1\}$.
The configuration space for the model is the set $\Lambda^N$
containing all
sequences $\omega= (\omega_1,\omega_2,
\ldots, \omega_N)$ with each $\omega_j \in\Lambda$. In terms of a
positive parameter $K$
representing the interaction strength, the Hamiltonian is defined by
\[
H_{N, K}(\omega)=\sum_{j=1}^N \omega_j^2 - \frac{K}{N} \Biggl( \sum_{j=1}^N
\omega_j \Biggr)^2
\]
for each $\omega\in\Lambda^N$.
Let $P_N$ be the product measure on $\Lambda^N$ with identical
one-dimensional marginals
$\rho= \frac{1}{3}(\delta_{-1}+\delta_0+\delta_1)$.
Thus $P_N$ assigns the probability $3^{-N}$ to each $\omega\in\Lambda^N$.
For inverse temperature $\beta> 0$ and for $K>0$,
the canonical ensemble for the mean-field B--C
model is the sequence of probability measures that assign to each
subset $B$ of $\Lambda^N$ the probability
%
\begin{eqnarray}
\label{eqn:pnbetak}
P_{N, \beta, K}(B) & = & \frac{1}{Z_N(\beta, K)} \cdot\int_B \exp
[-\beta H_{N,K}] \,dP_N \nonumber\\[-8pt]\\[-8pt]
& = & \frac{1}{Z_N(\beta, K)} \cdot\sum_{\omega\in B} \exp[-\beta
H_{N,K}(\omega)]
\cdot3^{-N}.
\nonumber
\end{eqnarray}
In this formula
$Z_N(\beta, K)$ is the partition function equal to
\[
\int_{\Lambda^N} \exp[-\beta H_{N,K}] \,dP_N
= \sum_{\omega\in\Lambda^N} \exp[-\beta H_{N,K}(\omega)] \cdot3^{-N}.
\]
Expectation with respect to $P_{N,\beta,K}$ is denoted by $E_{N,\beta,K}$.

The analysis of the
canonical ensemble $P_{N, \beta, K}$ is facilitated by absorbing
the noninteracting component of the Hamiltonian into the product
measure $P_N$, obtaining
%
\begin{equation}
\label{eqn:rewritecanon}
P_{N,\beta,K}(d\omega) =
\frac{1}{\tilde{Z}_N(\beta,K)} \cdot\exp\biggl[ N \beta K
\biggl(\frac{S_N(\omega)}{N} \biggr)^2 \biggr] P_{N,\beta}(d\omega).
\end{equation}
In this formula $S_N(\omega)$ equals
the total spin $\sum_{j=1}^N \omega_j$,
$P_{N,\beta}$ is the product measure on $\Lambda^N$ with identical
one-dimensional marginals
%
\begin{equation}
\label{eqn:rhobeta}
\rho_\beta(d \omega_j) =
\frac{1}{Z(\beta)} \cdot\exp(-\beta\omega_j^2) \rho(d \omega_j),
\end{equation}
$Z(\beta)$ is the normalization equal to
$\int_\Lambda\exp(-\beta\omega_j^2) \rho(d \omega_j) = (1 + 2
e^{-\beta})/3$
and $\tilde{Z}_N(\beta,K)$ is the normalization equal to
$Z_N(\beta,K)/[Z(\beta)]^N$.


We denote by $\mathcal{M}_{\beta,K}$ the set of equilibrium
macrostates of the
mean-field B--C model.
In order to describe this set, we
introduce the cumulant generating function $c_\beta$ of the measure
$\rho_\beta$ defined in (\ref{eqn:rhobeta}); for $t \in
{\mathbb R}
$ this function is defined by
\begin{eqnarray}
c_\beta(t) & = & \log\int_\Lambda\exp(t\omega_1) \rho_\beta
(d\omega_1)
\nonumber\\[-8pt]\\[-8pt]
& = &
\log\biggl( \frac{1+e^{-\beta}(e^t+e^{-t})}{1+2e^{-\beta}}
\biggr).\nonumber
\end{eqnarray}
For $x \in{\mathbb R}$ we define
%
\begin{equation}
\label{eqn:introgbetak}
G_{\beta,K}(x) = \beta K x^2 - c_\beta(2\beta K x).
\end{equation}
As shown in Proposition 3.4 in \cite{EllOttTou}, the set $\mathcal
{M}_{\beta,K}$
of equilibrium macrostates of the mean-field B--C model
can be characterized as the set of global minimum points of
$G_{\beta,K}$:
%
\begin{equation}
\label{eqn:ebetak}
\mathcal{M}_{\beta,K}= \{x \in[-1,1] \dvtx x \mbox{ is a global minimum
point of }
G_{\beta,K}(x)\}.
\end{equation}
In \cite{EllOttTou} the set $\mathcal{M}_{\beta,K}$
was denoted by $\tilde\mathcal{E}_{\beta,K}$.

We also define the canonical free energy
\[
\varphi(\beta,K) = - \lim_{N \rightarrow\infty} \frac{1}{\beta N}
\log
\tilde Z_N(\beta,K),
\]
where $\tilde Z_N(\beta,K)$ is the normalizing constant in
(\ref{eqn:rewritecanon}). This limit exists and
equals $\min_{x \in{\mathbb R}}\beta^{-1} G_{\beta,K}(x)$.
Because of this property of $G_{\beta,K}$, we call $G_{\beta,K}$ the
free-energy
functional of the mean-field B--C model.

The next two theorems use (\ref{eqn:ebetak}) to determine
the structure of $\mathcal{M}_{\beta,K}$ for $0 < \beta\leq\beta
_c= \log4$
and for $\beta> \beta_c$. The positive quantity $m(\beta,K)$ appearing
in these theorems
is called the thermodynamic magnetization.
The first theorem, proved in Theorem 3.6 in \cite{EllOttTou},
describes the continuous bifurcation in $\mathcal{M}_{\beta,K}$ for
$0 < \beta\leq
\beta_c$ as $K$ crosses a curve $\{(\beta,K(\beta))\dvtx 0 < \beta<
\beta_c\}$.
This bifurcation corresponds to a second-order phase transition, and
this curve is called the second-order curve.
The quantity $K(\beta)$, defined in (\ref{eqn:kcbeta}), is
denoted by $K^{(2)}_c(\beta)$ in \cite{EllOttTou}.
\begin{theorem}
\label{thm:secondorder}
For $0 < \beta\leq\beta_c$, we define
%
\begin{equation}
\label{eqn:kcbeta}
K(\beta) = {1}/[{2\beta c''_\beta(0)}] =
({e^\beta+ 2})/({4\beta}).
\end{equation}
For these values of $\beta$, $\mathcal{M}_{\beta,K}$ has the
following structure:

\textup{(a)} For $0 < K \leq K(\beta)$,
${\mathcal{M}}_{\beta,K} = \{0\}$.

\textup{(b)} For $K > K(\beta)$, there exists
${m}(\beta,K) > 0$ such that
${\mathcal{M}}_{\beta,K} = \{\pm m(\beta,K) \}$.

\textup{(c)} ${m}(\beta,K)$ is
a positive, increasing, continuous function for $K > K_c(\beta)$, and
as $K \rightarrow(K(\beta))^+$, $m(\beta,K) \rightarrow0$.
Therefore, ${\mathcal{M}}_{\beta,K}$ exhibits a continuous bifurcation
at $K(\beta)$.
\end{theorem}

The next theorem, proved in Theorem 3.8 in \cite{EllOttTou},
describes the discontinuous bifurcation in $\mathcal{M}_{\beta,K}$
for $\beta> \beta_c
$ as $K$ crosses
a curve $\{(\beta,K_1(\beta))\dvtx \beta> \beta_c\}$.
This bifurcation corresponds to a first-order phase transition, and
this curve
is called the first-order curve. As shown in Theorem 3.8 in
\cite{EllOttTou}, for all $\beta> \beta_c$, $K_1(\beta) < K(\beta)$. The
quantity $K_1(\beta)$ is denoted by $K^{(1)}_c(\beta)$ in \cite{EllOttTou}.
\begin{theorem}
\label{thm:firstorder}
For $\beta> \beta_c$, $\mathcal{M}_{\beta,K}$ has
the following structure in terms of the quantity
$K_1(\beta)$, denoted by $K_c^{(1)}(\beta)$ in \cite{EllOttTou}
and defined implicitly for $\beta> \beta_c$ on page 2231 of
\cite{EllOttTou}:

\textup{(a)} For $0 < K < K_1(\beta)$,
${\mathcal{M}}_{\beta,K} = \{0\}$.

\textup{(b)} For $K = K_1(\beta)$ there exists $m(\beta,K_1(\beta)) > 0$
such that ${\mathcal{M}}_{\beta,K_1(\beta)} =
\{0,\pm m(\beta,K_1(\beta))\}$.

\textup{(c)} For $K > K_1(\beta)$
there exists $m(\beta,K) > 0$
such that ${\mathcal{M}}_{\beta,K} =
\{\pm m(\beta,K)\}$.

\textup{(d)} $m(\beta,K)$
is a positive, increasing, continuous function for $K \geq K_1(\beta
)$, and
as $K \rightarrow K_1(\beta)^+$, $m(\beta,K) \rightarrow
m(\beta,K_1(\beta)) > 0$. Therefore,
${\mathcal{M}}_{\beta,K}$ exhibits a discontinuous bifurcation
at $K_1(\beta)$.
\end{theorem}

The phase-coexistence region is defined as the set of all points in the
positive $\beta$-$K$ quadrant for which
$\mathcal{M}_{\beta,K}$ consists of more than one value.
According to Theorems \ref{thm:secondorder} and \ref{thm:firstorder},
the phase-coexistence region consists of
all points above the second-order curve, above the tricritical point,
on the first-order curve
and above the first-order curve; that is,
\[
\{(\beta,K) \dvtx 0 < \beta\leq\beta_c, K > K(\beta) \mbox{ and }
\beta>
\beta_c, K \geq K_1(\beta)\}.
\]
Our derivation of the asymptotic behavior of the finite-size
magnetization\break
$E_{n,\beta_n,K_n}\{|S_n/n|\} \rightarrow0$ in this paper is valid
for a class of
sequences $(\beta_n,K_n)$ lying in the phase-coexistence region
for all sufficiently large $n$ and converging
either to a second-order point or to the tricritical point.
In the next section we state an asymptotic formula for $m(\beta_n,K_n)
\rightarrow0$
for a general class of such sequences.
That asymptotic formula will be used later in the paper when we study
the asymptotic behavior of the finite-size magnetization $E_{n,\beta
_n,K_n}\{
|S_n/n|\} \rightarrow0$.


\section{Asymptotic behavior of $m(\beta_n,K_n)$}
\label{section:mbnkn}

The main result in this section is Theorem
\ref{thm:exactasymptotics}. It states the asymptotic behavior of
the thermodynamic magnetization $m(\beta_n,K_n) \rightarrow0$
for sequences $(\beta_n,K_n)$ lying in the phase-coexistence region
for all sufficiently large $n$
and converging either to a second-order point or to the tricritical point.
The asymptotic behavior is expressed in terms of the unique positive,
global minimum point
of an associated polynomial that is introduced in hypothesis (iii) of
the theorem.
With several modifications the hypotheses
of the next theorem are also the hypotheses under which we derive the
rates at which
$E_{n,\beta_n,K_n}\{|S_n/n|\} \rightarrow0$ later in the paper.

As shown in part (iii) of Theorem \ref{thm:exactasymptotics}, the
asymptotics of $m(\beta_n,K_n)$ depend on the asymptotics
of the scaled free-energy function $n^{\alpha/\alpha_0}
G_{\beta_n,K_n}(x/n^{\theta\alpha})$. Because of Lemma \ref{lem:G},
the asymptotics of the finite-size magnetization in Theorems
\mbox{\ref{thm:mainsmall}--\ref{thm:mainmiddle}} depend on precisely
the same asymptotics. Lemma \ref{lem:G} coincides with Lem\-ma~4.1
in~\cite{CosEllOtt}. In that paper the connections
among the asymptotics of the scaled free-energy functional, the
limit theorems underlying the asymptotics
of the finite-size magnetization and Lemma 4.1 are described in detail.
These limit
theorems are analogues of the MDP in Theorem \ref{thm:mdp}
and of the weak convergence limit in Theorem \ref{thm:scalinglargealpha}.

Theorem \ref{thm:exactasymptotics} restates
the main theorem in \cite{EllMacOtt1}, Theorem 3.2.
Hypotheses (iii)(a) and (iv) in the next theorem coincide with
hypotheses (iii)(a) and (iv) in Theorem~3.2 in \cite{EllMacOtt1}
except that the latter hypotheses are expressed in terms of $u =
1-\alpha/\alpha_0$ and $\gamma= \theta\alpha$ while
here we have substituted the formulas for $u$ and~$\gamma$. Hence $u$
and $\gamma$ no longer appear.
%
%
\begin{theorem}
\label{thm:exactasymptotics}
Let $(\beta_n,K_n)$ be a positive sequence that converges either to a
second-order point
$(\beta,K(\beta))$, $0 < \beta< \beta_c$, or to the tricritical point
$(\beta,K(\beta)) = (\beta_c,K(\beta_c))$.
We assume that $(\beta_n,K_n)$ satisfies the following four
hypotheses:
\begin{longlist}
\item $(\beta_n,K_n)$ lies in the phase-coexistence region
for all sufficiently large $n$.

\item The sequence $(\beta_n,K_n)$ is parametrized by
$\alpha
> 0$. This parameter regulates the speed of approach of
$(\beta_n,K_n)$ to the second-order point or the tricritical point in the
following sense:
\[
b = \lim_{n \rightarrow\infty} n^\alpha(\beta_n- \beta)
\quad\mbox{and}\quad
k = \lim_{n \rightarrow\infty} n^\alpha\bigl(K_n- K(\beta)\bigr)
\]
both exist, and $b$ and $k$ are not both $0$; if $b \neq0$, then $b$
equals $1$ or $-1$.

\item There exists an even
polynomial $g$ of degree $4$ or $6$ satisfying $g(x) \rightarrow\infty
$ as
$|x| \rightarrow\infty$
together with the following two properties; $g$ is called the
Ginzburg--Landau polynomial.

\begin{longlist}[(a)\ \,(iii)]
\item[(a)] $\exists\alpha_0 > 0$ and $\exists\theta> 0$ such
that for all $\forall\alpha> 0$
\[
\lim_{n \rightarrow\infty} n^{\alpha/\alpha_0} G_{\beta
_n,K_n}(x/n^{\theta
\alpha}) = g(x)
\]
uniformly for $x$ in compact subsets of ${\mathbb R}$.

\item[(b)] $g$ has a unique, positive global minimum point $\bar{x}
$; thus the set of global minimum
points of $g$ equals $\{\pm\bar{x}\}$ or $\{0, \pm\bar{x}\}$.\vadjust{\goodbreak}
\end{longlist}

\item[(iv)] There exists a polynomial $H$ satisfying $H(x)
\rightarrow
\infty$ as $|x| \rightarrow\infty$
together with the following property: $\forall\alpha> 0$
$\exists R > 0$ such that
$\forall n \in\mathbb N$ sufficiently large and $\forall x \in
{\mathbb R}$
satisfying $|x/n^{\theta\alpha}| < R$, $n^{\alpha/\alpha_0}
G_{\beta_n,K_n}
(x/n^{\theta\alpha}) \geq H(x)$.
\end{longlist}
Under hypotheses \textup{(i)--(iv)}, for any $\alpha> 0$
\[
m(\beta_n,K_n)\sim{\bar{x}}/{n^{\theta\alpha}}\mbox{,\qquad that is, }
\lim_{n \rightarrow\infty} n^{\theta\alpha} m(\beta_n,K_n)= \bar{x}.
\]
If $b \not= 0$, then this becomes $m(\beta_n,K_n)\sim\bar{x}|\beta
- \beta_n
|^{\theta}$.
\end{theorem}

It is clear from the proof of the theorem that if hypotheses (iii) and
(iv) are valid for a specific value of $\alpha> 0$,
then we obtain the asymptotic formula $m(\beta_n,K_n)\sim{\bar
{x}}/{n^{\theta\alpha}}$ for that value of $\alpha$.

In the next section, we state the main results on the rates at which
$E_{n, \beta_n, K_n} \{|S_n/ n|\} \rightarrow0$ for small
$\alpha$ satisfying $0 < \alpha< \alpha_0$, for large $\alpha$
satisfying $\alpha> \alpha_0$,
and for intermediate $\alpha$ satisfying $\alpha= \alpha_0$.
We also compare these rates with the asymptotic behavior of the
thermodynamic magnetization $m(\beta_n,K_n)\rightarrow0$.


\section{Main results on rates at which $E_{n, \beta_n, K_n} \{
|S_n/n|\} \rightarrow0$}
\label{section:rates}

Let $a_n$ be a positive sequence converging to 0. In stating the three
results on the rates at which
the finite-size magnetization $E_{n, \beta_n, K_n} \{|S_n/n|\}
\rightarrow
0$, we write
\[
E_{n, \beta_n, K_n}\{|S_n/n|\} \sim a_n \qquad\mbox{if }
\lim_{n \rightarrow\infty} E_{n, \beta_n, K_n}\{|S_n/n|\}/a_n = 1,
\]
and we write
\[
E_{n, \beta_n, K_n}\{|S_n/n|\} \gg a_n \qquad\mbox{if }
\lim_{n \rightarrow\infty} E_{n, \beta_n, K_n}\{|S_n/n|\}/a_n =
\infty.
\]

Let $\alpha$ be the quantity parametrizing the sequences $(\beta_n,K_n)$
as explained in hypothesis (ii)
of Theorem \ref{thm:exactasymptotics}.
We begin with Theorem \ref{thm:mainsmall}, which gives the rate at
which $E_{n, \beta_n, K_n}\{|S_n/n|\} \rightarrow0$ for small $\alpha$
satisfying $0 < \alpha< \alpha_0$.
Theorem \ref{thm:mainlarge}
gives the rate at which $E_{n, \beta_n, K_n}\{|S_n/n|\} \rightarrow0$ for
large $\alpha$ satisfying $\alpha> \alpha_0$ while
Theorem \ref{thm:mainmiddle} gives the rate at which $E_{n, \beta_n,
K_n}\{|S_n/n|\} \rightarrow0$
for intermediate $\alpha$ satisfying $\alpha= \alpha_0$. In all
three cases we compare these rates with the rate
at which $m(\beta_n,K_n)\rightarrow0$.
In the next section we specialize these theorems to the six sequences
mentioned in the \hyperref[section:intro]{Introduction}.

Part (a) of the next theorem gives the rate at which $E_{n, \beta_n,
K_n}\{|S_n/n|\} \rightarrow0$
for $0 < \alpha< \alpha_0$, and part (b) shows that for these values
of $\alpha$,
$E_{n, \beta_n, K_n}\{|S_n/n|\} \sim m(\beta_n,K_n)$. It follows that
for $0
< \alpha< \alpha_0$, $m(\beta_n,K_n)$ is a
physically relevant estimator of the finite-size magnetization $E_{n,
\beta_n, K_n}\{|S_n/n|\}$
because it has the same asymptotic behavior as that quantity.

The next theorem is valid under hypotheses (i) and (ii) of Theorem \ref
{thm:exactasymptotics}, hypotheses (iii)(a) and (iv) of
that theorem for all $0 < \alpha< \alpha_0$, the inequality $0 <
\theta\alpha_0 < 1/2$,
and a new hypothesis (iii$^\prime$)(b). The inequality $0 < \theta
\alpha_0 < 1/2$ is satisfied by all six
sequences considered in Section \ref{section:sequences}.
The new hypothesis (iii$^\prime$)(b) restricts hypothesis (iii)(b) of
Theorem \ref{thm:exactasymptotics}
by assuming that the set of global minimum points of the
Ginzburg--Landau polynomial $g$ equals $\{\pm\bar{x}\}$
for some $\bar{x}$. As we remark after the statement of the theorem,
this restriction is needed in order
to prove part (a). The proof does not cover the case where the set of
global minimum points of $g$ equals $\{0, \pm\bar{x}\}$ for some
$\bar{x}
> 0$. The conjecture is that in this case
there exists $0 < \lambda< 1/2$ such that $E_{n,\beta_n,K_n}\{
|S_n/n|\} \sim2\lambda\bar{x}/n^{\theta\alpha}$ (see the discussion
before Corollary \ref{cor:weaklimit}).
An example of a sequence for which the set of
global minimum points of $g$ contains three points is given in case (d)
of sequence 4 in the next section.
By contrast, all the other sequences considered in the next section satisfy
the new hypothesis that the set of global minimum points of $g$ equals
$\{\pm\bar{x}\}$
for some $\bar{x}$.
\begin{theorem}[($0 < \alpha< \alpha_0$)]
\label{thm:mainsmall}
Let $(\beta_n,K_n)$ be a positive sequence para\-metrized by $\alpha> 0$
and converging either to a second-order point $(\beta,K(\beta))$,
$0 < \beta< \beta_c$, or to the tricritical point $(\beta_c,K(\beta
_c))$. We
assume hypotheses \textup{(i)} and \textup{(ii)}
of Theorem \ref{thm:exactasymptotics} together with hypotheses
\textup{(iii)(a)} and \textup{(iv)}
of that theorem for all $0 < \alpha< \alpha_0$. We also assume the
inequality $0 < \theta\alpha_0 < 1/2$
and the following hypothesis, which restricts
hypothesis \textup{(iii)(b)} of Theorem~\ref{thm:exactasymptotics}:
\begin{longlist}[(iii$^\prime$)(b)]
\item[(iii$^\prime$)(b)] The set of global minimum points of the
Ginzburg--Landau polynomial $g$ equals $\{\pm\bar{x}\}$
for some $\bar{x}> 0$.
\end{longlist}
The following conclusions hold:

\textup{(a)} For all $0 < \alpha< \alpha_0$
\[
E_{n,\beta_n, K_n}\{|S_n/n|\} \sim{\bar{x}}/{n^{{\theta\alpha}}}\mbox{,\qquad
that is, } \lim_{n \rightarrow\infty} n^{\theta\alpha} E_{n,
\beta
_n, K_n} \{ | S_n/n | \} = \bar{x}.
\]

\textup{(b)} For all $0 < \alpha< \alpha_0$, $E_{n,\beta_n, K_n}\{
|S_n/n|\} \sim m(\beta_n, K_n)$.
\end{theorem}

Part (a) of the theorem is proved from the moderate deviation
principle (MDP) for the $P_{n,\beta_n,K_n}$-distributions of
$S_n/n^{1-\theta
\alpha}$ in Theorem \ref{thm:mdp}, which shows that the rate
function equals $g - \inf_{y \in{\mathbb R}}g(y)$. The inequality $0
< \theta
\alpha_0 < 1/2$ is used to control an error
term in the proof of the MDP. According to hypothesis (iii$^\prime
$)(b), the set
of global minimum points of $g$ equals $\{\pm\bar{x}\}$ for some
$\bar{x}
> 0$. It quickly follows from the MDP
that the sequence of $P_{n,\beta_n,K_n}$-distributions of
$S_n/n^{1-\theta
\alpha}$ converges
weakly to $\frac{1}{2}\delta_{\bar{x}} + \frac{1}{2}\delta_{-\bar
{x}}$. The
uniform integrability of $S_n/n^{1-\theta\alpha}$,
derived in Lemma~\ref{lem:uniformint} from the MDP, yields the limit
$E_{n,\beta_n,K_n}\{|S_n/n^{1-\theta\alpha}|\} \rightarrow\bar{x}$
as $n
\rightarrow\infty$.
This is the asymptotic
formula for $E_{n,\beta_n,K_n}\{|S_n/n|\}$ in part (a) of Theorem \ref
{thm:mainsmall}.
Part (b) of the theorem follows from part (a) and the asymptotic formula
$m(\beta_n,K_n) \sim\bar{x}/n^{\theta\alpha}$, which is the
conclusion of Theorem \ref{thm:exactasymptotics}.


We next state Theorem \ref{thm:mainlarge}, which in part (a) gives
the rate at which\break $E_{n, \beta_n, K_n}\{|S_n/n|\} \rightarrow0$ for
$\alpha
> \alpha_0$.
Part (b) shows that for these values of $\alpha$,
$E_{n, \beta_n, K_n}\{|S_n/n|\} \gg m(\beta_n,K_n)$.
Because $m(\beta_n,K_n)\rightarrow0$ at an asymptotically faster rate
than the finite-size magnetization $E_{n,\beta_n,K_n}\{|S_n/n|\}$,
$m(\beta_n,K_n)$ is not a physically relevant estimator of that
quantity for
$\alpha> \alpha_0$.

In order to prove part (a) of the next theorem, we need hypothesis (iv)
of Theorem \ref{thm:exactasymptotics} for $\alpha= \alpha_0$,
the inequality $0 < \theta\alpha_0 < 1/2$ and a new hypothesis (v),
in which we assume
that for all $\alpha> \alpha_0$, $nG_{\beta_n,K_n}(x/n^{\theta
\alpha_0})$
convergence pointwise to a polynomial $\tilde{g}(x)$ that goes to
$\infty$
as $|x| \rightarrow\infty$.
As we will see for the first five of the six sequences considered in
the next section,
$\tilde{g}$ in hypothesis (v) equals the highest order term of the
Ginzburg--Landau polynomial $g$.
We omit the analysis showing that this description of $\tilde{g}$ can,
in fact, be validated in general
if the uniform convergence in hypothesis (iii)(a) of Theorem \ref
{thm:exactasymptotics} on compact subsets of ${\mathbb R}$
is strengthened to uniform convergence on compact subsets of an
appropriate open set in $\mathbb C$ containing the origin and if
$\theta
\alpha_0$ equals a certain value depending on the degree
of $g$. This stronger convergence is valid for the six sequences considered
in the next section. However, the additional condition on $\theta
\alpha_0$, valid for the
first five sequences, is not satisfied by the sixth sequence.

In part (b) of the next theorem the rates at which
$E_{n, \beta_n, K_n}\{|S_n/n|\} \rightarrow0$ and $m(\beta
_n,K_n)\rightarrow0$
are compared. In order to prove part (b), we also need hypotheses (i)
and (ii) of Theorem \ref{thm:exactasymptotics}
and hypotheses (iii) and (iv) of that theorem for all $\alpha> \alpha
_0$. These hypotheses allow us to
apply Theorem \ref{thm:exactasymptotics} for all $\alpha> \alpha_0$.
%
%
\begin{theorem}[($\alpha> \alpha_0$)]
\label{thm:mainlarge}
Let $(\beta_n,K_n)$ be a positive sequence parametrized by $\alpha> 0$
and converging either to a second-order point $(\beta,K(\beta))$,
$0 < \beta< \beta_c$, or to the tricritical point $(\beta_c,K(\beta
_c))$. We
assume hypotheses \textup{(i)} and \textup{(ii)}
of Theorem~\ref{thm:exactasymptotics}, hypothesis \textup{(iii)} of
Theorem \ref{thm:exactasymptotics} for
all $\alpha> \alpha_0$ and hypothesis \textup{(iv)} of Theorem
\ref{thm:exactasymptotics} for all $\alpha\geq\alpha_0$. We
also assume the inequality
$0 < \theta\alpha_0 < 1/2$ and the following hypothesis:
\begin{longlist}
\item[(v)] There exists an even polynomial $\tilde{g}$ of degree
$4$ or $6$ satisfying
$\tilde{g}(x) \rightarrow\infty$ as $|x| \rightarrow\infty$
together with the
following property:
$\exists\alpha_0 > 0$ and $\exists\theta> 0$ such that $\forall
\alpha> \alpha_0$ and $\forall x \in{\mathbb R}$
\[
\lim_{n \rightarrow\infty} n G_{\beta_n,K_n}(x/n^{\theta\alpha
_0}) = \tilde{g}(x).
\]
\end{longlist}
The following conclusions hold:

\textup{(a)} We define
\[
\bar{y} = \frac{1}{\int_{{\mathbb R}} \exp[-\tilde{g}(x)] \,dx}
\cdot\int
_{{\mathbb R}} |x| \exp[-\tilde{g}(x)] \,dx.
\]
Then for all $\alpha> \alpha_0$
\[
E_{n, \beta_n, K_n}\{|S_n/n|\} \sim{\bar{y}}/{n^{\theta\alpha_0}}\mbox{,\qquad
that is, }
\lim_{n \rightarrow\infty} n^{\theta\alpha_0} E_{n, \beta_n,
K_n}\{
|S_n/n|\} = \bar{y}.
\]

\textup{(b)} For all $\alpha> \alpha_0$, $E_{n, \beta_n, K_n}\{
|S_n/n|\} \gg m(\beta_n,K_n)$.
\end{theorem}

Part (a) of the theorem is proved from the weak convergence
of the sequence of $P_{n,\beta_n,K_n}$-distributions of
$S_n/n^{1-\theta\alpha
}$ to a probability measure having a density
proportional to $\exp[-\tilde{g}]$, which is shown
in Theorem \ref{thm:scalinglargealpha}. The proof of this weak
convergence relies on hypothesis (v) of Theorem \ref{thm:mainlarge}
and the lower bound in hypothesis (iv) of Theorem \ref
{thm:exactasymptotics} for $\alpha= \alpha_0$.
The inequality $0 < \theta\alpha_0 < 1/2$ is used to control an error term
in the proof. The uniform-integrability-type estimate in Proposition
\ref{prop:weakerunifint} yields the limit $E_{n,\beta_n,K_n}\{
|S_n/n^{1-\theta\alpha_0}|\} \rightarrow\bar{y}$
as $n \rightarrow\infty$. This is~the asymptotic
formula for $E_{n,\beta_n,K_n}\{|S_n/n|\}$ in part (a) of Theorem \ref
{thm:mainlarge}.
Part (b) of the theorem follows from part (a), the asymptotic
formula $m(\beta_n,K_n) \sim\bar{x}/n^{\theta\alpha}$ and the fact
that since $\alpha> \alpha_0$,
the decay rate $n^{-\theta\alpha}$ of $m(\beta_n,K_n)\rightarrow0$
is asymptotically larger than the decay rate $n^{-\theta\alpha_0}$ of
$E_{n, \beta_n, K_n}\{|S_n/n|\}$.

We end this section by stating Theorem \ref{thm:mainmiddle}.
Part (a) gives the rate at which $E_{n, \beta_n, K_n}\{|S_n/n|\}
\rightarrow0$
for $\alpha= \alpha_0$, and part (b) compares this rate
with the rate at which $m(\beta_n,K_n) \rightarrow0$.
The theorem is valid under hypotheses (i) and (ii) of Theorem \ref
{thm:exactasymptotics},
hypotheses (iii) and (iv) of Theorem \ref{thm:exactasymptotics} for
$\alpha= \alpha_0$
and the inequality $0 < \theta\alpha_0 < 1/2$.
%
%
\begin{theorem}[($\alpha= \alpha_0$)]
\label{thm:mainmiddle}
Let $(\beta_n,K_n)$ be a positive sequence parametrized by $\alpha> 0$
and converging either to a second-order point $(\beta,K(\beta))$,
$0 < \beta< \beta_c$, or to the tricritical point $(\beta_c,K(\beta
_c))$. We
assume hypotheses \textup{(i)} and \textup{(ii)}
of Theorem \ref{thm:exactasymptotics}, hypotheses \textup{(iii)}
and \textup{(iv)}
of Theorem \ref{thm:exactasymptotics} for $\alpha= \alpha_0$
and the inequality $0 < \theta\alpha_0 < 1/2$.
The following conclusions hold:

\textup{(a)} We define
\[
\bar{z} = \frac{1}{\int_{{\mathbb R}} \exp[-g(x)] \,dx} \cdot\int
_{{\mathbb R}} |x|
\exp[-g(x)] \,dx.
\]
Then for all $\alpha= \alpha_0$
\[
E_{n, \beta_n, K_n}\{|S_n/n|\} \sim{\bar{z}}/{n^{\theta\alpha_0}}\mbox{,
\qquad that is, }
\lim_{n \rightarrow\infty} n^{\theta\alpha_0} E_{n, \beta_n,
K_n}\{
|S_n/n|\} = \bar{z}.
\]

\textup{(b)} For $\alpha= \alpha_0$, $E_{n, \beta_n, K_n}\{|S_n/n|\}
\sim\bar{z} \cdot m(\beta_n,K_n)/\bar{x}$.
\end{theorem}

We omit the proof of part (a) of the theorem, which can be derived like
part (a) of Theorem
\ref{thm:mainlarge}. According to hypothesis (iii)(a) of Theorem
\ref{thm:exactasymptotics}, $n G_{\beta_n,K_n}(x/n^{\theta\alpha_0})$
converges to $g(x)$
uniformly for $x$ in compact subsets of ${\mathbb R}$. The pointwise
convergence of $n G_{\beta_n,K_n}(x/n^{\theta\alpha_0})$ to $g(x)$
and the lower bound in hypothesis (iv) of Theorem \ref{thm:exactasymptotics}
for $\alpha= \alpha_0$
allow us to prove that the sequence of $P_{n,\beta
_n,K_n}$-distributions of
$S_n/n^{1-\theta\alpha_0}$
converges weakly to a probability measure having a density
proportional to $\exp[-{g}]$. The inequality $0 < \theta\alpha_0 <
1/2$ is used to control an error term
in the proof. The asymptotic
formula for $E_{n,\beta_n,K_n}\{|S_n/n|\}$ in part (a) of Theorem~\ref
{thm:mainmiddle} follows from this weak-convergence limit
and the uniform-integrability-type estimate in Proposition \ref
{prop:weakerunifint}, the hypotheses of which can be verified
in the context of Theorem \ref{thm:mainmiddle} as they are verified at
the end of Section \ref{section:largealpha} in the context
of Theorem \ref{thm:mainlarge}. When $\alpha= \alpha_0$,
$m(\beta_n,K_n) \sim\bar{x}/n^{\theta\alpha_0}$ [Theorem \ref
{thm:exactasymptotics}(b)]. Hence part (a) of Theorem
\ref{thm:mainmiddle} implies that
\[
E_{n, \beta_n, K_n}\{|S_n/n|\} \sim{\bar{z}}/{n^{\theta\alpha_0}}
\sim\bar{z} \cdot m(\beta_n,K_n)/\bar{x}.
\]
This is the conclusion of part (b) of the theorem.

In numerical calculations we studied the relative size of $\bar{z}$ and
$\bar{x}$.
Depending on the magnitude of the coefficient of the quadratic term in
the Ginzburg--Landau polynomial $g$,
$\bar{z}/\bar{x}$ can be less than 1, can equal 1 and can exceed 1.


In the next section we specialize Theorems \ref{thm:mainsmall}, \ref
{thm:mainlarge} and \ref{thm:mainmiddle}
to the six sequences mentioned in the \hyperref[section:intro]{Introduction}.


\section{Results for six sequences}
\label{section:sequences}

In \cite{EllMacOtt1} we apply Theorem \ref{thm:exactasymptotics} to determine
the asymptotic behavior of the thermodynamic magnetization
$m(\beta_n,K_n)\rightarrow0$ for six sequences $(\beta_n,K_n)$
parametrized by $\alpha> 0$. The first two sequences converge to a
second-order point $(\beta,K(\beta))$,
$0 < \beta< \beta_c$, and the last four sequences converge to the
tricritical point $(\beta_c,K(\beta_c))$.
In the present section we specialize to the first five sequences the results
in Theorems \ref{thm:mainsmall}, \ref{thm:mainlarge} and \ref
{thm:mainmiddle} concerning the
the asymptotic behaviors of $E_{n, \beta_n, K_n}\{|S_n/n|\}
\rightarrow0$
for $0 < \alpha< \alpha_0$, $\alpha> \alpha_0$
and $\alpha= \alpha_0$. We also compare these asymptotic behaviors
with the
asymptotic behavior of $m(\beta_n,K_n) \rightarrow0$. In addition we state
the results of Theorems \ref{thm:exactasymptotics}, \ref{thm:mainsmall}
and~\ref{thm:mainmiddle} for the sixth sequence. However, for this
sequence, one of the hypotheses of Theorem \ref{thm:mainlarge}
is not valid, and so that theorem cannot be applied.

In order to be able to apply these four theorems, we must verify the
validity of their hypotheses, which are the following:
\begin{itemize}
\item \textit{Theorem} \ref{thm:exactasymptotics}. Hypotheses (i) and (ii)
and hypotheses (iii) and (iv) for all $\alpha> 0$.
\item \textit{Theorem} \ref{thm:mainsmall}. Hypotheses (i) and (ii) of
Theorem \ref{thm:exactasymptotics}, hypotheses (iii)(a) and (iv)
of Theorem \ref{thm:exactasymptotics} for all $0 < \alpha< \alpha
_0$, the inequality
$0 < \theta\alpha_0 < 1/2$ and the new hypothesis (iii$^\prime$)(b).
\item \textit{Theorem} \ref{thm:mainlarge}. Hypotheses (i) and (ii) of
Theorem \ref{thm:exactasymptotics},
hypothesis (iii) of Theorem \ref{thm:exactasymptotics} for all $\alpha
> \alpha_0$, hypothesis
(iv) of Theorem \ref{thm:exactasymptotics} for all $\alpha\geq\alpha
_0$, the inequality
$0 < \theta\alpha_0 < 1/2$ and the new hypothesis (v).
\item \textit{Theorem} \ref{thm:mainmiddle}. Hypotheses (i) and (ii) of
Theorem \ref{thm:exactasymptotics},
hypotheses (iii) and (iv) of Theorem \ref{thm:exactasymptotics}
for $\alpha= \alpha_0$ and the inequality $0 < \theta\alpha_0 < 1/2$.
\end{itemize}
Thus, in order to verify the hypotheses of the four theorems, it
suffices to verify hypotheses (i) and (ii) of Theorem
\ref{thm:exactasymptotics}, hypotheses (iii)(a) and (iv) of Theorem~\ref
{thm:exactasymptotics}
for all $\alpha> 0$, hypothesis (iii$^\prime$)(b) for all $0 < \alpha < \alpha_0$,
hypothesis (iii)(b) for all $\alpha \geq \alpha_0$, the inequality $0 < \theta\alpha_0 < 1/2$, and
hypothesis (v) of Theorem \ref{thm:mainlarge}.

The quantities $\bar{y}$ and $\bar{z}$ appearing in the
asymptotic formulas in Theorems \ref{thm:mainlarge} and
\ref{thm:mainmiddle} are defined as follows in terms of the polynomial
$\tilde{g}$,
introduced in hypothesis (v) of Theorem \ref{thm:mainlarge}, and in
terms of the Ginzburg--Landau polynomial $g$:
%
\[
\bar{y}= \frac{1}{\int_{{\mathbb R}} \exp[-\tilde{g}(x)] \,dx} \cdot
\int_{{\mathbb R}
} |x| \exp[-\tilde{g}(x)] \,dx
\]
and
\[
\bar{z} = \frac{1}{\int_{{\mathbb R}} \exp[-g(x)] \,dx} \cdot\int
_{{\mathbb R}} |x|
\exp[-g(x)] \,dx.
\]
For the first five sequences, $\tilde{g}$ equals the highest-order
term in $g$. For the sixth sequence,
Theorem \ref{thm:mainlarge} cannot be applied because hypothesis (v)
of that theorem is not valid.
In each sequence $K(\beta) = (e^\beta+ 2)/(4\beta)$ for
$\beta> 0$. The curve $\{(\beta,K(\beta))\dvtx 0 < \beta< \beta_c\}$ is
the second-order curve, $(\beta_c,K(\beta_c))$ is the
tricritical point and the curve $\{(\beta,K(\beta))\dvtx \beta> \beta
_c\}$
is the spinodal curve.\vspace*{8pt}

\textit{Sequence} 1.

\textit{Definition of sequence} 1. Given $0 < \beta< \beta_c$, $\alpha> 0$,
$b \in\{1,0,-1\}$
and $k \in{\mathbb R}$, $k \not= 0$, the sequence is defined by
%
\begin{equation}
\label{eqn:bnone} \beta_n = \beta+ b/n^\alpha\quad\mbox{and}\quad K_n =
K(\beta) + k/n^\alpha.
\end{equation}
This sequence converges to the second-order point $(\beta,K(\beta))$
along a ray with slope $k/b$ if $b \not= 0$.\vspace*{8pt}

\textit{Hypotheses} (i) \textit{and} (ii) \textit{in Theorem} \ref{thm:exactasymptotics}.
Hypothesis (i) states that $(\beta_n,K_n)$ lies in the phase-coexistence region
for all sufficiently large $n$. In order to guarantee this, we assume that
$K'(\beta)b - k < 0$. This inequality
is equivalent to $K_n > K(\beta_n)$ for all sufficiently large $n$ and
thus guarantees that
$(\beta_n,K_n)$ lies in the phase-coexistence region above the
second-order curve for all sufficiently large $n$. Hypothesis (ii) is
also satisfied.\vspace*{8pt}

\textit{Other hypotheses}.
\begin{enumerate}
\item Define $\alpha_0 = 1/2$ and $\theta= 1/2$. As shown in Theorem
4.1 in \cite{EllMacOtt1},
the uniform convergence in hypothesis (iii)(a) of Theorem
\ref{thm:exactasymptotics} is valid for all $\alpha> 0$ with the
Ginzburg--Landau polynomial
\begin{eqnarray}
g(x) = \beta\bigl(K'(\beta) b - k\bigr) x^2 + c_4(\beta) x^4\nonumber\\
\eqntext{\mbox{where }
c_4(\beta) = (e^\beta+ 2)^2(4 - e^\beta)/8\cdot4!.}
\end{eqnarray}
Since $\theta\alpha_0 = 1/4$, we have $0 < \theta\alpha_0 < 1/2$,
which is one of the hypotheses
of Theorems \ref{thm:mainsmall}--\ref{thm:mainmiddle}.
\item We assume that $K'(\beta)b - k < 0$.
Then, as required by hypothesis (iii)(b) of Theorem \ref
{thm:exactasymptotics} and hypothesis (iii$^\prime$)(b) of Theorem
\ref{thm:mainsmall}, the set of global minimum points of $g$ is $\{\pm
\bar{x}\}$, where $\bar{x}> 0$
is defined in (4.6) in \cite{EllMacOtt1}.
\item Hypothesis (iv) of Theorem \ref{thm:exactasymptotics} is valid
for all $\alpha> 0$ with the polynomial $H$ given on page 113 of
\cite{EllMacOtt1}.
\item The pointwise convergence in hypothesis (v) of Theorem \ref
{thm:mainlarge}
holds with $\tilde{g}$ equal to the highest order term in $g$; namely,
$\tilde{g}(x) = c_4(\beta) x^4$. This is easily verified using
equation (4.4) in \cite{EllMacOtt1}.
\end{enumerate}

We now specialize to sequence 1 the results in
Theorems \ref{thm:exactasymptotics}, \ref{thm:mainsmall}, \ref
{thm:mainlarge} and \ref{thm:mainmiddle} concerning the
asymptotic behavior of $m(\beta_n,K_n) \rightarrow0$ and the
asymptotic behaviors of $E_{n, \beta_n, K_n}\{|S_n/n|\} \rightarrow0$ for
$0 < \alpha< \alpha_0$, for $\alpha> \alpha_0$
and for $\alpha= \alpha_0$.
\begin{theorem}
\label{thm:resultsone}
Let $(\beta_n, K_n)$ be sequence 1 that is defined in
(\ref{eqn:bnone}) and converges
to a second-order point $(\beta,K(\beta))$ for $0 < \beta< \beta_c$.
Assume that
$K'(\beta)b - k < 0$. The following conclusions hold:

\textup{(a)} For all $\alpha> 0$,
\[
m(\beta_n,K_n) \sim\bar{x}/n^{\alpha/2}.
\]
If $b \not= 0$ in the definition of $\beta_n$, then $m(\beta_n,K_n)
\sim
\bar{x}|\beta- \beta_n|^{1/2}$.

\textup{(b)} For all $0 < \alpha< \alpha_0 = 1/2$,
\[
E_{n, \beta_n, K_n}\{| S_n/n|\} \sim{\bar{x}}/{n^{\alpha/2}} \sim
m(\beta_n, K_n).
\]

\textup{(c)} For all $\alpha> \alpha_0 = 1/2$,
\[
E_{n, \beta_n, K_n}\{|S_n/n|\} \sim{\bar{y}}/{n^{1/4}} \gg m(\beta
_n, K_n).
\]

\textup{(d)} For $\alpha= \alpha_0 = 1/2$,
\[
E_{n, \beta_n, K_n}\{| S_n/n|\} \sim{\bar{z}}/{n^{1/4}} \sim\bar{z}
\cdot m(\beta_n, K_n)/\bar{x}.
\]
\end{theorem}

\textit{Sequence} 2.

\textit{Definition of sequence} 2. Given $0 < \beta_0 < \beta_c$,
$\alpha> 0$, $b \in\{1,-1\}$, an integer $p \geq2$ and a real number
$\ell\not= K^{(p)}(\beta)$, the sequence is defined by
%
\begin{eqnarray}
\label{eqn:bntwo}
\beta_n&=& \beta_0 + {b}/{n^\alpha} \quad\mbox{and}\nonumber\\[-8pt]\\[-8pt]
K_n&=& K(\beta_0) + \sum_{j=1}^{p-1} {K^{(j)}(\beta_0) b^j}/(j!
n^{j\alpha}) +
{\ell b^p}/(p! n^{p\alpha}).\nonumber
\end{eqnarray}
This sequence converges to the second-order point $(\beta_0,K(\beta_0))$
along a curve that coincides with the second-order curve
to order $p-1$ in powers of $\beta- \beta_0$.\vspace*{8pt}

\textit{Hypotheses} (i) \textit{and} (ii) \textit{in Theorem} \ref{thm:exactasymptotics}.
Hypothesis (i)
states that $(\beta_n,K_n)$ lies in the phase-coexistence region
for all sufficiently large $n$. In order to guarantee this, we assume that
$(K^{(p)}(\beta_0) - \ell)b^p < 0$. This inequality
is equivalent to $K_n > K(\beta_n)$ for all sufficiently large $n$ and
thus guarantees that
$(\beta_n,K_n)$ lies in the phase-coexistence region above the
second-order curve for all sufficiently large $n$.
Hypothesis (ii) is also satisfied.\vspace*{8pt}

\textit{Other hypotheses}.
\begin{enumerate}
\item Define $\alpha_0 = 1/2p$ and $\theta= p/2$.
As shown in Theorem 4.2 in \cite{EllMacOtt1},
the uniform convergence in hypothesis (iii)(a) of Theorem
\ref{thm:exactasymptotics} is valid for all $\alpha> 0$ with the
Ginzburg--Landau polynomial
\begin{eqnarray}
g(x) &=& \frac{1}{p!}\beta_0\bigl(K^{(p)}(\beta_0) - \ell\bigr)b^p x^2 +
c_4(\beta_0) x^4\nonumber\\
\eqntext{\mbox{where } c_4(\beta_0) = (e^{\beta_0} + 2)^2(4 - e^{\beta
_0})/8\cdot4!.}
\end{eqnarray}
Since $\theta\alpha_0 = 1/4$, we have $0 < \theta\alpha_0 < 1/2$,
which is one of the hypotheses of
Theorems \ref{thm:mainsmall}--\ref{thm:mainmiddle}.
\item We assume that $(K^{(p)}(\beta_0) - \ell)b^p < 0$.
Then, as required by hypothesis (iii)(b) of Theorem \ref
{thm:exactasymptotics} and hypothesis (iii$^\prime$)(b) of
Theorem \ref{thm:mainsmall},
the set of global minimum points of $g$ is $\{\pm\bar{x}\}$, where
$\bar{x}> 0$ is defined in (4.9) in \cite{EllMacOtt1}.
\item Hypothesis (iv) of Theorem \ref{thm:exactasymptotics} is valid
for all $\alpha> 0$ with the polynomial $H$ given on page 115 of
\cite{EllMacOtt1}.
\item The pointwise convergence in hypothesis (v) of Theorem \ref
{thm:mainlarge}
holds with $\tilde{g}$ equal to the highest order term in $g$; namely,
$\tilde{g}(x) = c_4(\beta_0) x^4$. This is easily verified using
(4.8) in \cite{EllMacOtt1}.
\end{enumerate}

We now specialize to sequence 2 the results in
Theorems \ref{thm:exactasymptotics}, \ref{thm:mainsmall}, \ref
{thm:mainlarge} and \ref{thm:mainmiddle} concerning the
asymptotic behavior of $m(\beta_n,K_n) \rightarrow0$ and the
asymptotic behaviors $E_{n, \beta_n, K_n}\{|S_n/n|\} \rightarrow0$
for $0 <
\alpha< \alpha_0$, for $\alpha> \alpha_0$
and for $\alpha= \alpha_0$.
\begin{theorem}
\label{thm:resultstwo}
Let $(\beta_n, K_n)$ be sequence 2 that is defined in
(\ref{eqn:bntwo}) and converges
to a second-order point $(\beta_0,K(\beta_0))$ for $0 < \beta_0 <
\beta_c$. Assume that
$(K^{(p)}(\beta_0) - \ell)b^p < 0$. The following conclusions hold:

\textup{(a)} For all $\alpha> 0$,
\[
m(\beta_n,K_n) \sim\bar{x}/n^{p\alpha/2} = \bar{x}|\beta_0 -
\beta_n|^{p/2}.
\]

\textup{(b)} For all $0 < \alpha< \alpha_0 = 1/2p$,
\[
E_{n, \beta_n, K_n}\{| S_n/n|\} \sim{\bar{x}}/{n^{p\alpha/2}} \sim
m(\beta_n, K_n).
\]

\textup{(c)} For all $\alpha> \alpha_0 = 1/2p$,
\[
E_{n, \beta_n, K_n}\{|S_n/n|\} \sim{\bar{y}}/{n^{1/4}} \gg m(\beta
_n, K_n).
\]

\textup{(d)} For $\alpha= \alpha_0 = 1/2p$,
\[
E_{n, \beta_n, K_n}\{| S_n/n|\} \sim{\bar{z}}/{n^{1/4}} \sim\bar{z}
\cdot m(\beta_n, K_n)/\bar{x}.
\]
\end{theorem}

\textit{Sequence} 3.

\textit{Definition of sequence} 3. This sequence is defined as in
(\ref{eqn:bnone}) with $\beta$ replaced by $\beta_c$.
Thus given $\alpha> 0$, $b \in\{1,0,-1\}$,
and $k \in{\mathbb R}$, $k \not= 0$, the sequence is defined by
%
\begin{equation}
\label{eqn:bnthree}
\beta_n = \beta_c+ b/n^\alpha\quad\mbox{and}\quad K_n = K(\beta_c) +
k/n^\alpha.
\end{equation}
This sequence converges to the tricritical point $(\beta_c,K(\beta
_c))$ along
a ray with slope $k/b$ if $b \not= 0$.\vspace*{8pt}

\textit{Hypotheses} (i) \textit{and} (ii) \textit{in Theorem} \ref
{thm:exactasymptotics}.
Hypothesis (i) states that $(\beta_n,K_n)$ lies in the phase-coexistence region
for all sufficiently large $n$. In order to guarantee this, we assume that
$K'(\beta_c)b - k < 0$. This inequality
is equivalent to $K_n > K(\beta_n)$ for all sufficiently large $n$ and
thus guarantees that for all sufficiently large $n$,
$(\beta_n,K_n)$ lies in the phase-coexistence region above the spinodal
curve if $b = 1$, above the second-order curve
if $b = -1$ and above the tricritical point if $b = 0$. Hypothesis~(ii)
is also satisfied.\vadjust{\goodbreak}

\textit{Other hypotheses}.
\begin{enumerate}
\item Define $\alpha_0 = 2/3$ and $\theta= 1/4$. As shown in Theorem
4.3 in \cite{EllMacOtt1},
the uniform convergence in hypothesis (iii)(a) of Theorem
\ref{thm:exactasymptotics} is valid for all $\alpha> 0$ with the
Ginzburg--Landau polynomial
\[
g(x) = \beta_c\bigl(K'(\beta_c) b - k\bigr) x^2 + c_6 x^6\qquad \mbox{where } c_6
= 9/40.
\]
Since $\theta\alpha_0 = 1/6$, we have $0 < \theta\alpha_0 < 1/2$,
which is one of the hypotheses
of Theorems \ref{thm:mainsmall}--\ref{thm:mainmiddle}.
\item We assume that $K'(\beta_c)b - k < 0$. Then, as required by
hypothesis (iii)(b)
of Theorem \ref{thm:exactasymptotics} and hypothesis (iii$^\prime
$)(b) of Theorem \ref{thm:mainsmall},
the set of global minimum points of $g$ is $\{\pm\bar{x}\}$, where
$\bar{x}> 0$ is defined in (4.14) in \cite{EllMacOtt1}.
\item Hypothesis (iv) of Theorem \ref{thm:exactasymptotics} is valid
for all $\alpha> 0$ with the polynomial $H$ given on page 117 of
\cite{EllMacOtt1}.
\item The pointwise convergence in hypothesis (v) of Theorem \ref
{thm:mainlarge}
holds with $\tilde{g}$ equal to the highest order term in $g$; namely,
$\tilde{g}(x) = c_6 x^6$.
This is easily verified using (4.13) in \cite{EllMacOtt1}.
\end{enumerate}

We now specialize to sequence 3 the results in
Theorems \ref{thm:exactasymptotics}, \ref{thm:mainsmall}, \ref
{thm:mainlarge} and \ref{thm:mainmiddle} concerning the
asymptotic behavior of $m(\beta_n,K_n) \rightarrow0$ and the
asymptotic behaviors of $E_{n, \beta_n, K_n}\{|S_n/n|\} \rightarrow0$ for
$0 < \alpha< \alpha_0$, for $\alpha> \alpha_0$
and for $\alpha= \alpha_0$.
\begin{theorem}
\label{thm:resultsthree}
Let $(\beta_n, K_n)$ be sequence 3 that is defined in
(\ref{eqn:bnthree}) and
converges to the tricritical point $(\beta_c,K(\beta_c))$. Assume that
$K'(\beta_c)b - k < 0$. The following conclusions hold:

\textup{(a)} For all $\alpha> 0$,
\[
m(\beta_n,K_n) \sim\bar{x}/n^{\alpha/4}.
\]
If $b \not= 0$ in the definition of $\beta_n$, then $m(\beta_n,K_n)
\sim
\bar{x}|\beta- \beta_n|^{1/4}$.

\textup{(b)} For all $0 < \alpha< \alpha_0 = 2/3$,
\[
E_{n, \beta_n, K_n}\{| S_n/n|\} \sim{\bar{x}}/{n^{\alpha/4}} \sim
m(\beta_n, K_n).
\]

\textup{(c)} For all $\alpha> \alpha_0 = 2/3$,
\[
E_{n, \beta_n, K_n}\{|S_n/n|\} \sim{\bar{y}}/{n^{1/6}} \gg m(\beta
_n, K_n).
\]

\textup{(d)} For $\alpha= \alpha_0 = 2/3$,
\[
E_{n, \beta_n, K_n}\{| S_n/n|\} \sim{\bar{z}}/{n^{1/6}} \sim\bar{z}
\cdot m(\beta_n, K_n)/\bar{x}.
\]
\end{theorem}

\textit{Sequence} 4.

Of the six sequences this sequence exhibits the most complicated
behavior, the description of which is divided into four cases (a)--(d)
described in the third paragraph below.
In addition, for cases (c) and (d) the validity of hypothesis (i) of
Theorem \ref{thm:exactasymptotics} involves the validity
of two conjectures. For cases (a), (b) and (c) all the other hypotheses
of Theorems
\ref{thm:exactasymptotics}, \ref{thm:mainsmall}, \ref{thm:mainlarge}
and \ref{thm:mainmiddle} are valid.
However, for case (d) hypothesis (iii$^\prime$)(b) of Theorem \ref
{thm:mainsmall} is not valid,
and therefore that theorem cannot be applied in that case.\vspace*{8pt}

\textit{Definition of sequence} 4. Given $\alpha> 0$, a curvature
parameter $\ell\in{\mathbb R}$, and another
parameter $\tilde\ell\in{\mathbb R}$, sequence 4 is defined by
%
\begin{equation}
\label{eqn:bnfour}
\beta_n = \beta_c+ {1}/{n^\alpha} \quad\mbox{and}\quad K_n = K(\beta_c) +
K'(\beta_c
)/{n^\alpha}
+ \ell/(2n^\alpha) + \tilde\ell/(6n^{3\alpha}).\hspace*{-33pt}
\end{equation}
Since $\beta_n- \beta_c= 1/n^\alpha$ the sequence converges from the right
to the tricritical point $(\beta_c,K(\beta_c))$
along the curve $(\beta,\tilde{K}(\beta))$, where for $\beta> \beta_c$
\[
\tilde{K}(\beta) = K(\beta_c)+ K'(\beta_c)(\beta- \beta_c) + \ell
(\beta-
\beta_c)^2/2 +
\tilde\ell(\beta- \beta_c)^3/6.
\]
At the tricritical point this curve is tangent to the spinodal curve,
which is the extension of the second-order
curve to $\beta> \beta_c$. As shown in \cite{EllOttTou}, Theorem 3.8,
the spinodal curve lies above the first-order curve for all $\beta>
\beta_c$.\vspace*{8pt}

\textit{Hypotheses} (i) \textit{and} (ii) \textit{in Theorem}
\ref{thm:exactasymptotics}. The discussion of hypothesis (i) for this
sequence involves four cases (a), (b), (c) and (d) that are presented
in the next paragraph. The validity of this hypotheses for the last two
of these four cases depends on the validity of conjectures 1 and 2
stated at the end of this paragraph. These conjectures are supported by
partial proofs, numerical evidence and properties of the
Ginzburg--Landau polynomials and are discussed in detail in Section 6
of \cite{EllMacOtt2}. The two conjectures involve the behavior, in a
neighborhood of the tricritical point, of the first-order curve defined
by $K_1(\beta)$ for $\beta> \beta_c$. Since
$\lim_{\beta\rightarrow\beta_c^+} K_1(\beta) = K(\beta_c)$~\cite{EllOttTou},
Sections 3.1 and 3.3, by continuity we extend the
definition of $K_1(\beta)$ to $\beta_c$ by defining $K_1(\beta _c) =
K(\beta_c)$. We assume that the first three right-hand derivatives of
$K_1(\beta)$ exist at $\beta_c$ and denote them by $K_1'(\beta_c)$,
$K_1''(\beta_c)$ and $K_1'''(\beta_c)$. We also define $\ell_c =
K''(\beta_c) - {5}/({4\beta_c})$. Conjectures 1 and 2 state the
following: (1) $K_1'(\beta_c) = K'(\beta_c)$, (2) $K_1''(\beta_c) =
\ell_c < 0 < K''(\beta_c)$.

The choices of $\ell$ and $\tilde\ell$ defining the four cases of sequence
4 are as follows.
Cases (a)--(c) correspond to $\ell> \ell_c$ and suitable values of
$\tilde\ell$,
and case (d) corresponds to $\ell= \ell_c$ and suitable values of
$\tilde\ell$.
%
\begin{itemize}
\item[(a)] $\ell> K''(\beta_c)$ and any $\tilde\ell\in{\mathbb R}$.
\item[(b)] $\ell= K''(\beta_c)$ and any $\tilde\ell> K'''(\beta_c)$.
\item[(c)] $K''(\beta_c) > \ell> \ell_c$ and any $\tilde\ell\in
{\mathbb R}$.
\item[(d)] $\ell= \ell_c$ and any $\tilde\ell> K_1'''(\beta_c)$.
\end{itemize}

For all four cases hypothesis (ii) is satisfied. For cases (a) and (b)
and for
all sufficiently large $n$, $(\beta_n,K_n)$ lies in the phase-coexistence
region above the spinodal curve,
and so hypothesis (i) is satisfied.
If conjectures 1 and 2 are valid, then for cases (c) and (d) $(\beta_n,K_n
)$ lies in the phase-coexistence region between the
spinodal and first-order curves for all sufficiently large $n$, and
again hypothesis (i) is valid.
In the discussion of the validity of the hypotheses
of Theorems \ref{thm:exactasymptotics} and \ref{thm:mainsmall}--\ref
{thm:mainmiddle} for sequence 4,
conjectures 1 and 2 are needed only for the last assertion. For all
four cases
$(\beta_n,K_n)$ converges to the tricritical point along the curve
$\{(\beta,\tilde{K}(\beta))\dvtx \beta> \beta_c\}$, where $\tilde
{K}(\beta)$
is defined in the display after (\ref{eqn:bnfour}).
If conjectures 1 and 2 are valid, then for cases (a)--(c) this curve
coincides with the first-order curve to order 1 in powers of $\beta-
\beta_c$,
while for case (d) this curve coincides with the first-order curve to
order 2 in powers of $\beta- \beta_c$.\vspace*{8pt}

\textit{Other hypotheses.}

The validity of these hypotheses for cases (c)
and (d) does not
depend on conjectures 1 and 2. A major difference between cases
(a)--(c) and case (d) appears in item 2.
\begin{enumerate}
\item Define $\alpha_0 = 1/3$ and $\theta= 1/2$. As shown in Theorem
4.4 in \cite{EllMacOtt1},
for all four cases the uniform convergence in hypothesis (iii)(a) of Theorem
\ref{thm:exactasymptotics} is valid for all $\alpha> 0$ with the
Ginzburg--Landau polynomial
\begin{eqnarray}
g(x) = \tfrac{1}{2}\beta_c\bigl(K^{\prime\prime}(\beta_c) - \ell\bigr)x^2 - 4 c_4 x^4 + c_6
x^6\nonumber\\
\eqntext{\mbox{where } c_4 = 3/16 \mbox{ and } c_6 = 9/40.}
\end{eqnarray}
Since $\theta\alpha_0 = 1/6$, we have $0 < \theta\alpha_0 < 1/2$,
which is one of the hypotheses of
Theorems \ref{thm:mainsmall}--\ref{thm:mainmiddle}.
\item We assume that $\ell> \ell_c = K''(\beta_c) - 5/(4\beta_c)$.
Then, as required by hypothesis~(iii)(b) of Theorem \ref{thm:exactasymptotics}
and hypothesis (iii$^\prime$)(b)
of Theorem \ref{thm:mainsmall}, for cases (a)--(c)
the set of global minimum points of $g$ equals $\{\pm\bar{x}\}$, where
$\bar{x}= \bar{x}(\ell) > 0$ is defined in (4.19) in \cite{EllMacOtt1}.
If $\ell= \ell_c$, then for case (d) the set of global minimum points of
$g$ equals $\{0,\pm\bar{x}\}$, where $\bar{x}= \bar{x}(\ell_c) >
0$ is
defined in (4.19) in~\cite{EllMacOtt1}.
Hence for case (d) hypothesis (iii)(b) of Theorem \ref
{thm:exactasymptotics} is valid, but hypothesis (iii$^\prime$)(b)
of Theorem \ref{thm:mainsmall} is not valid.
\item For all four cases,
hypothesis (iv) of Theorem \ref{thm:exactasymptotics} is valid for all
$\alpha> 0$ with the polynomial $H$ given on page 120 of
\cite{EllMacOtt1}.
\item For all four cases, the pointwise convergence in hypothesis (v)
of Theorem \ref{thm:mainlarge}
holds with $\tilde{g}$ equal to the highest order term in $g$; namely,
$\tilde{g}(x) = c_6 x^6$.
This is easily verified using (4.16) in \cite{EllMacOtt1}.
\end{enumerate}

We now specialize to sequence 4 the results in
Theorems \ref{thm:exactasymptotics}, \ref{thm:mainsmall}, \ref
{thm:mainlarge} and \ref{thm:mainmiddle} concerning the
asymptotic behavior of $m(\beta_n,K_n) \rightarrow0$ and the
asymptotic behaviors of $E_{n, \beta_n, K_n}\{|S_n/n|\} \rightarrow0$ for
$0 < \alpha< \alpha_0$, for $\alpha> \alpha_0$
and for $\alpha= \alpha_0$. Parts (a), (c) and (d) of the theorem are
valid for all four cases of the sequence. However, part (b)
is valid only for cases (a), (b) and (c) because, as we point out in
item 2 above,
for case (d) hypothesis (iii$^\prime$)(b) of Theorem \ref
{thm:mainsmall} does not hold.
\begin{theorem}
\label{thm:resultsfour}
Let $(\beta_n, K_n)$ be sequence 4 that is defined in
(\ref{eqn:bnfour}) and converges
to the tricritical point $(\beta_c,K(\beta_c))$.
Assume that $\ell$ and $\tilde\ell$ are defined as in one of the
four cases
\textup{(a)--(d)} and that for cases \textup{(c)--(d)} conjectures 1
and 2 are valid. The following conclusions hold:

\textup{(a)} For cases \textup{(a)--(d)}, for all $\alpha> 0$,
\[
m(\beta_n,K_n) \sim\bar{x}/n^{\alpha/2} = \bar{x}(\beta_n - \beta
_c)^{1/2}.
\]

\textup{(b)} For cases \textup{(a)--(c)}, for all $0 < \alpha< \alpha_0 =
1/3$,
\[
E_{n, \beta_n, K_n}\{| S_n/n|\} \sim{\bar{x}}/{n^{\alpha/2}} \sim
m(\beta_n, K_n).
\]

\textup{(c)} For cases \textup{(a)--(d)}, for all $\alpha> \alpha_0 =
1/3$,
\[
E_{n, \beta_n, K_n}\{|S_n/n|\} \sim{\bar{y}}/{n^{1/6}} \gg m(\beta
_n, K_n).
\]

\textup{(d)} For cases \textup{(a)--(d)}, for $\alpha= \alpha_0 = 1/3$,
\[
E_{n, \beta_n, K_n}\{| S_n/n|\} \sim{\bar{z}}/{n^{1/6}} \sim\bar{z}
\cdot m(\beta_n, K_n)/\bar{x}.
\]
\end{theorem}

\textit{Sequence} 5.

\textit{Definition of sequence} 5. This sequence is defined as in
(\ref{eqn:bntwo})
with $b = -1$, $p = 2$ and $\beta_0$ replaced by $\beta_c$. Thus given
$\alpha> 0$ and a real number $\ell\not= K^{\prime\prime}(\beta_c)$,
the sequence is defined by
%
\begin{equation}
\label{eqn:bnfive}
\beta_n= \beta_c- 1/{n^\alpha} \quad\mbox{and}\quad
K_n= K(\beta_c) - K^{\prime}(\beta_c)/n^\alpha+ \ell/2n^{2\alpha}.
\end{equation}
This sequence converges to the tricritical point $(\beta_c,K(\beta
_c))$ from
the left
along a curve that coincides with the second-order curve to order $2$
in powers of $\beta- \beta_c$.\vspace*{8pt}

\textit{Hypotheses} (i) \textit{and} (ii) \textit{in Theorem} \ref{thm:exactasymptotics}.
Hypothesis (i) states that $(\beta_n,K_n)$ lies in the phase-coexistence region
for all sufficiently large $n$. In order to guarantee this, we assume that
$\ell> K^{\prime\prime}(\beta_c)$. This inequality
is equivalent to $K_n > K(\beta_n)$ for all sufficiently large $n$ and
thus guarantees that
$(\beta_n,K_n)$ lies in the phase-coexistence region above the
second-order curve for all sufficiently large $n$.
Hypothesis (ii) is also satisfied.\vspace*{8pt}

\textit{Other hypotheses}.
\begin{enumerate}
\item Define $\alpha_0 = 1/3$ and $\theta= 1/2$.
As shown in Theorem 4.5 in \cite{EllMacOtt1},
the uniform convergence in hypothesis (iii)(a) of Theorem
\ref{thm:exactasymptotics} is valid for all $\alpha> 0$ with the
Ginzburg--Landau polynomial
\begin{eqnarray}
g(x) = \tfrac{1}{2}\beta_c\bigl(K^{\prime\prime}(\beta_c) - \ell\bigr)x^2 + 4 c_4 x^4 + c_6
x^6\nonumber\\
\eqntext{\mbox{where } c_4 = 3/16 \mbox{ and } c_6 = 9/40.}
\end{eqnarray}
Since $\theta\alpha_0 = 1/6$, we have $0 < \theta\alpha_0 < 1/2$,
which is one of the hypotheses of
Theorems \ref{thm:mainsmall}--\ref{thm:mainmiddle}.
\item We assume that $\ell> K''(\beta_c)$. Then, as required by
hypothesis (iii)(b)
of Theorem \ref{thm:exactasymptotics} and hypothesis (iii$^\prime
$)(b) of
Theorem \ref{thm:mainsmall},
the set of global minimum points of $g$ is $\{\pm\bar{x}\}$, where
$\bar{x}> 0$ is defined in (4.23) in \cite{EllMacOtt1}.
\item Hypothesis (iv) of Theorem \ref{thm:exactasymptotics} is valid
for all $\alpha> 0$ with the polynomial $H$ given on page 121 of
\cite{EllMacOtt1}.
\item The pointwise convergence in hypothesis (v) of Theorem \ref
{thm:mainlarge}
holds with $\tilde{g}$ equal to the highest order term in $g$; namely,
$\tilde{g}(x) = c_6 x^6$.
This is easily verified using (4.21) in \cite{EllMacOtt1}.
\end{enumerate}

We now specialize to sequence 5 the results in
Theorems \ref{thm:exactasymptotics}, \ref{thm:mainsmall}, \ref
{thm:mainlarge} and \ref{thm:mainmiddle} concerning the
asymptotic behavior of $m(\beta_n,K_n) \rightarrow0$ and the
asymptotic behaviors of $E_{n, \beta_n, K_n}\{|S_n/n|\} \rightarrow0$ for
$0 < \alpha< \alpha_0$, for $\alpha> \alpha_0$
and for $\alpha= \alpha_0$.
\begin{theorem}
\label{thm:resultsfive}
Let $(\beta_n, K_n)$ be sequence 5 that is defined in
(\ref{eqn:bnfive}) and converges
to the tricritical point $(\beta_c,K(\beta_c))$. Assume that $\ell>
K''(\beta_c)$.
The following conclusions hold:

\textup{(a)} For all $\alpha> 0$,
\[
m(\beta_n,K_n) \sim\bar{x}/n^{\alpha/2} = \bar{x}(\beta_c- \beta
_n)^{1/2}.
\]

\textup{(b)} For all $0 < \alpha< \alpha_0 = 1/3$,
\[
E_{n, \beta_n, K_n}\{| S_n/n|\} \sim{\bar{x}}/{n^{\alpha/2}} \sim
m(\beta_n, K_n).
\]

\textup{(c)} For all $\alpha> \alpha_0 = 1/3$,
\[
E_{n, \beta_n, K_n}\{|S_n/n|\} \sim{\bar{y}}/{n^{1/6}} \gg m(\beta
_n, K_n).
\]

\textup{(d)} For $\alpha= \alpha_0 = 1/3$,
\[
E_{n, \beta_n, K_n}\{| S_n/n|\} \sim{\bar{z}}/{n^{1/6}} \sim\bar{z}
\cdot m(\beta_n, K_n)/\bar{x}.
\]
\end{theorem}

\textit{Sequence} 6.

For this sequence the hypotheses of Theorems \ref{thm:exactasymptotics},
\ref{thm:mainsmall} and \ref{thm:mainmiddle} are all valid.
However, Theorem \ref{thm:mainlarge} cannot be applied because
hypothesis (v) of that theorem is not valid.\vspace*{8pt}

\textit{Definition of sequence} 6. This sequence is defined as in
(\ref{eqn:bntwo})
with $b = -1$, an integer $p \geq3$, and $\beta_0$ replaced by $\beta_c
$. Thus given
$\alpha> 0$ and a real number $\ell\not= K^{(p)}(\beta_c)$, the
sequence is defined by
%
\begin{eqnarray}
\label{eqn:bnsix}
\beta_n&=& \beta_c - 1/{n^\alpha} \quad\mbox{and}\nonumber\\[-8pt]\\[-8pt]
K_n&=& K(\beta_c) + \sum_{j=1}^{p-1} {K^{(j)}(\beta_c) (-1)^j}/(j!
n^{j\alpha
}) +
{\ell(-1)^p}/(p! n^{p\alpha}).\nonumber
\end{eqnarray}
This sequence converges to the tricritical point $(\beta_c,K(\beta
_c))$ from
the left
along a curve that coincides with the second-order curve to order $p-1$
in powers of $\beta- \beta_c$.\vspace*{8pt}

\textit{Hypotheses} (i) \textit{and} (ii) \textit{in Theorem} \ref{thm:exactasymptotics}.
Hypothesis (i) states that $(\beta_n,K_n)$ lies in the phase-coexistence region
for all sufficiently large $n$. In order to guarantee this, we assume that
$(K^{(p)}(\beta_c) - \ell)(-1)^p < 0$. This inequality
is equivalent to $K_n > K(\beta_n)$ for all sufficiently large $n$ and
thus guarantees that
$(\beta_n,K_n)$ lies in the phase-coexistence region above the
second-order curve for all sufficiently large~$n$.
Hypothesis (ii) is also satisfied.\vspace*{8pt}

\textit{Other hypotheses}.
\begin{enumerate}
\item Define $\alpha_0 = 1/(2p-1)$ and $\theta= (p-1)/2$.
As shown in Theorem 4.6 in \cite{EllMacOtt1},
the uniform convergence in hypothesis (iii)(a) of Theorem
\ref{thm:exactasymptotics} is valid for all $\alpha> 0$ with the
Ginzburg--Landau polynomial
\[
g(x) = \frac{1}{p!}\beta_c\bigl(K^{(p)}(\beta_c) - \ell\bigr)(-1)^p x^2 + 4 c_4
x^4 \qquad \mbox{where } c_4 = 3/16.
\]
Since $\theta\alpha_0 = (p-1)/[2(2p-1]$ and $p \geq3$,
we have $0 < \theta\alpha_0 < 1/2$, which is one of the hypotheses of
Theorems \ref{thm:mainsmall}--\ref{thm:mainmiddle}.
\item We assume that $(K^{(p)}(\beta_c) - \ell)(-1)^p < 0$.
Then, as required by hypothesis (iii)(b) of Theorem
\ref{thm:exactasymptotics} and hypothesis (iii$^\prime$)(b) of
Theorem \ref{thm:mainsmall},
the set of global minimum points of $g$ is $\{\pm\bar{x}\}$, where
$\bar{x}> 0$ is defined in (4.25) in \cite{EllMacOtt1}.
\item Hypothesis (iv) of Theorem \ref{thm:exactasymptotics} is valid
for all $\alpha> 0$ with the polynomial $H$ given on page 122 of
\cite{EllMacOtt1}.
\item The only problem arises in hypothesis (v) of Theorem \ref{thm:mainlarge},
which is not valid for all $\alpha> \alpha_0$
with the values of $\alpha_0$ and $\theta$ in item 1. In fact, one uses
equation (4.21) in \cite{EllMacOtt1} to verify that with these values of
$\alpha_0$, $\theta$ and $\alpha$,
$nG_{\beta_n,K_n}(x/n^{\theta\alpha_0}) \rightarrow0$ for all $x
\in{\mathbb R}$.
Hence with these values of $\alpha_0$, $\theta$ and $\alpha$,
Theorem \ref{thm:mainlarge} cannot be applied.
\end{enumerate}

We now specialize to sequence 6 the results in
Theorems \ref{thm:exactasymptotics}, \ref{thm:mainsmall} and \ref
{thm:mainmiddle} concerning the
asymptotic behavior of $m(\beta_n,K_n) \rightarrow0$ and the
asymptotic behaviors
of $E_{n, \beta_n, K_n}\{|S_n/n|\} \rightarrow0$ for $0 < \alpha<
\alpha
_0$ and for $\alpha= \alpha_0$.
\begin{theorem}
\label{thm:resultssix}
Let $(\beta_n, K_n)$ be sequence 6 that is defined in
(\ref{eqn:bnfive}) and converges
to the tricritical point $(\beta_c,K(\beta_c))$. Assume that
$(K^{(p)}(\beta
) - \ell)(-1)^p < 0$.
The following conclusions hold:

\textup{(a)} For all $\alpha> 0$,
\[
m(\beta_n,K_n) \sim\bar{x}/n^{(p-1)\alpha/2} = \bar{x}(\beta_c-
\beta_n)^{(p-1)/2}.
\]

\textup{(b)} For all $0 < \alpha< \alpha_0 = 1/(2p-1)$,
\[
E_{n, \beta_n, K_n}\{| S_n/n|\} \sim{\bar{x}}/{n^{(p-1)\alpha/2}}
\sim m(\beta_n, K_n).
\]

\textup{(c)} For $\alpha= \alpha_0 = 1/(2p-1)$,
\[
E_{n, \beta_n, K_n}\{| S_n/n|\} \sim{\bar{z}}/{n^{(p-1)/[2(2p-1)]}}
\sim\bar{z}\cdot m(\beta_n, K_n)/\bar{x}.
\]
\end{theorem}

The one gap in Theorem \ref{thm:resultssix} is the failure of
hypothesis (v) of
Theorem \ref{thm:mainlarge} for all $\alpha> \alpha_0$. We omit the analysis
that gives a variation of Theorem \ref{thm:mainlarge} describing a subset
of $\alpha> \alpha_0$ for which the asymptotics of $E_{n,\beta
_n,K_n}\{
|S_n/n|\} \rightarrow0$
can be determined.

This completes our description, in the context of the six sequences,
of the three theorems in Section \ref{section:rates} on how
the asymptotic behaviors of the thermodynamic magnetization $m(\beta_n,K_n)
\rightarrow0$ and
the finite-size magnetization $E_{n,\beta_n,K_n}\{|S_n/n|\}
\rightarrow0$
compare for $0 < \alpha< \alpha_0$,
$\alpha> \alpha_0$ and $\alpha= \alpha_0$.
In the next section we outline the theory of finite-size scaling,
which gives insight into the physical phenomena underlying the theorems
in Section \ref{section:rates}.

\section{The theory of finite-size scaling}
\label{section:fss}

\quad In Theorems \ref{thm:mainsmall} and \ref{thm:mainlarge} we compare
the asymptotic behavior of the
thermodynamic magnetization\break $m(\beta_n,K_n)\rightarrow0$ with the asymptotic
behavior of the finite-size
magnetization $E_{n,\beta_n,K_n}\{|S_n/n|\} \rightarrow0$, first for
$0 < \alpha<
\alpha_0$ and then for $\alpha> \alpha_0$.
The results described in these two theorems
are intimately connected with the theory of finite-size scaling. This
nonrigorous but highly suggestive theory
was developed in statistical mechanics in order to understand phase
transitions in finite systems.
In fact, our work in this paper was motivated by the theory of
finite-size scaling
and can be understood in that context. At the same time, our results
put ideas of
finite-size scaling on a firm mathematical footing for the mean-field
B--C model.
To the best of our knowledge, this is the first time that the theory
of finite-size scaling has been rigorously derived for a mean-field
model. After sketching the theory of finite-size scaling, we show that its
predictions are consistent with those in Theorem \ref{thm:resultsone}.
That theorem specializes Theorems \ref{thm:mainsmall}
and \ref{thm:mainlarge} to sequence 1, which is defined in~(\ref{eqn:bnone}).

The theory of finite-size scaling is a generalization of scaling theory
to apply to finite systems \cite{Barber}.
Scaling theory gives a methodology for analyzing the singularities of
thermodynamic quantities such as the magnetization
in a neighborhood of criticality. One formulation of scaling theory
emphasizes the fundamental role of the correlation length $\xi$ by
expressing the singularities in thermodynamic quantities in terms of
$\xi$. For example, in a neighborhood
of criticality the thermodynamic magnetization behaves like $\xi
^{-\tilde{\beta}/\nu}$,
where $\tilde{\beta}$ is the magnetization exponent and $\nu$ is the
correlation-length exponent \cite{Stanley71}. The singularity in the
correlation length as a function of the distance to criticality
is controlled by the exponent $-\nu$.

The theory of finite-size scaling asserts that in a neighborhood of
criticality
quantities such as the finite-size magnetization behave like functions
of the linear system size $L$
and the ratio of the correlation length $\xi$ to the linear system
size. When $\xi/L \ll1$, the system is effectively infinite so that
finite-size quantities are independent of $L$, and the critical
singularities are the same as those in the thermodynamic limit.
On the other hand, when $\xi/L \gg1$, critical fluctuations are
instead limited by the system size. In this regime, the theory of
finite-size scaling asserts that the power-law singularities as a
function of $\xi$ are replaced by power-law
singularities as a function of $L$. For example, in the case of the
finite-size magnetization
the theory of finite-size scaling asserts that in a neighborhood of criticality
it behaves like $L^{-\tilde{\beta}/\nu} f(\xi/L)$. The function $f(x)$
interpolates continuously between the two regimes.
Thus, as $x = \xi/L \rightarrow0$, $f(x) \approx x^{-\tilde{\beta
}/\nu}$. In
this case the finite-size magnetization
behaves like $L^{-\tilde{\beta}/\nu}(\xi/L)^{-\tilde{\beta}/\nu
} = \xi^{-\tilde{\beta}
/\nu}$ and so is independent of $L$.
As discussed in the preceding paragraph, the thermodynamic
magnetization also behaves like the same function $\xi^{-\tilde
{\beta}/\nu}$.
On the other hand, as $x = \xi/L \rightarrow\infty$, $f(x)
\rightarrow
1$ and the finite-size magnetization behaves like
$L^{-\tilde{\beta}/\nu}$.

These ideas cannot be directly applied to the mean-field B--C model or
other mean-field
spin systems since neither the system length $L$ nor the correlation
length $\xi$ are defined. Appropriate quantities for mean-field
spin systems are $N$, the number of spins, and $\Xi$, the size of the
giant cluster in the Fortuin--Kasteleyn
representation \cite{ForKas,Gri}. For such systems the mappings
$N = L^{d_c}$ and $\Xi= \xi^{d_c}$ are expected to yield,
in a neighborhood of criticality, correct scaling relations for
thermodynamic quantities such as the magnetization and correct finite-size
scaling relations for quantities such as the finite-size magnetization.
In these equations $d_c$ denotes the upper critical dimension.
This is defined as the dimension above which
short-range spin systems such as the B--C model
\cite{Blu,Cap1,Cap2,Cap3}
have the same critical exponents as the associated mean-field models.
Thus in the case of the thermodynamic magnetization
the scaling expression $\xi^{-\tilde{\beta}/\nu}$, which is
appropriate for
short-range models, is replaced by $\Xi^{-\tilde{\beta}/d_c \nu}$.
In addition, in the case of the finite-size magnetization, the
finite-size scaling expression $L^{-\tilde{\beta}/\nu} f(\xi/L)$, which
is appropriate for short-range models, is replaced by $N^{-\tilde
{\beta}/d_c
\nu} f((\Xi/N)^{1/d_c})$.

In order to apply the ideas of finite-size scaling to the mean-field
B--C model, we consider a sequence
$(\beta_n,K_n)$ converging to criticality---that is, a second-order
point or the tricritical point---from the phase-coexistence region.
We also identify the number of spins $N$ with the index $n$
parametrizing the sequence $(\beta_n,K_n)$.
Thus the finite-size scaling expression for the finite-size
magnetization takes the form
$n^{-\tilde{\beta}/d_c \nu} f((\Xi/n)^{1/d_c})$. As in Section 5
of \cite{EllMacOtt1},
we bring in the quantity $\mu_1(\beta_n,K_n)$ representing the distance
of $(\beta_n,K_n)$ to criticality. According to scaling theory,
$\Xi$ behaves like $\mu_1^{-d_c\nu}$.

We now specialize these ideas to sequence 1.
Defined in (\ref{eqn:bnone}), this sequence
converges to a second-order point and $\mu_1\approx n^{-\alpha}$. Thus
for this sequence the correlation volume $\Xi$ behaves like
$\mu_1^{-d_c\nu} = n^{d_c\alpha\nu}$, and so the ratio $\Xi/n$
appearing in the argument\vspace*{1pt} of $f$ behaves like $n^{d_c\alpha\nu-1}$.
Since for mean-field second-order points $\nu=\tilde{\beta}=1/2$ and $d_c=4$
\cite{PlBe}, we see that $\Xi$ and $\Xi/n$
behave, respectively, like $n^{2\alpha}$ and $n^{2 \alpha-1}$. The
conclusion is that for sequence 1
the scaling relation for the thermodynamic magnetization takes the form
%
\begin{equation}
\label{eqn:scalingmagn}
\Xi^{-\tilde{\beta}/d_c \nu} \approx n^{-\alpha/2},
\end{equation}
and the finite-size scaling expression for the finite-size
magnetization takes the form
%
\begin{equation}
\label{eqn:fss}
n^{-\tilde{\beta}/d_c \nu} f\bigl((\Xi/n)^{1/d_c}\bigr)\approx n^{-1/4}
f\bigl(n^{(2 \alpha-1)/4}\bigr).
\end{equation}

The next step is to relate this phenomenology with the conclusions of
Theorem~\ref{thm:resultsone}, which
specializes Theorems \ref{thm:mainsmall} and \ref{thm:mainlarge} to
sequence 1. The key is to recall that
$f(x) \approx x^{-\tilde{\beta}/\nu} = x^{-1}$ as $x \rightarrow0$
and $f(x) \rightarrow
1$ as $x \rightarrow\infty$.
According to the formula in (\ref{eqn:fss}), the theory of
finite-size scaling predicts a change in behavior
in the finite-size magnetization when $\alpha= 1/2$.
This agrees with Theorem \ref{thm:resultsone}, which states that for sequence
1 the threshold value $\alpha_0$ equals $1/2$. For $0 < \alpha< 1/2$,
the ratio $\Xi/n = n^{2\alpha- 1}$
is much less than 1, and the finite-size magnetization behaves like
$n^{-1/4} n^{-(2 \alpha-1)/4} = n^{-\alpha/2}$.
This behavior coincides with the behavior of the thermodynamic
magnetization given in (\ref{eqn:scalingmagn}),
making this prediction
of the theory of finite-size scaling consistent with part (b) of
Theorem \ref{thm:resultsone}.
On the other hand, for $\alpha> 1/2$, since the ratio $\Xi/n =
n^{2\alpha-1}$ is much bigger than 1,
we have $f(n^{(2\alpha-1)/4}) \approx1$, and so the finite-size magnetization
behaves like $n^{-1/4}$. This converges to 0 much more slowly than the
thermodynamic magnetization,
which behaves like $n^{-\alpha/2}$. Again this prediction of the
theory of finite-size scaling is consistent
with part (c) of Theorem \ref{thm:resultsone}.

Similar heuristic arguments based on the theory of finite-size scaling
can be applied to the other sequences discussed in Section \ref
{section:sequences}.
They yield the correct asymptotic behaviors for the finite-size
magnetization for $0 < \alpha< \alpha_0$
and $\alpha> \alpha_0$, in agreement with Theorems \ref
{thm:resultstwo}--\ref{thm:resultssix}.
However, the tricritical region presents additional difficulties
because of the cross-over from the second-order
regime to the tricritical regime. The correct treatment of these
sequences in the scaling regime
is discussed in the context of scaling theory in Section 5 of
\cite{EllMacOtt1}.

This completes our discussion of the theory of finite-size scaling and
its relationship
with the main mathematical results given in Theorems \ref
{thm:mainsmall} and \ref{thm:mainlarge}
and specialized to the six sequences in Theorems \ref
{thm:resultsone}--\ref{thm:resultssix}.
In the next section we discuss how part (a) of Theorem \ref
{thm:mainsmall} follows from the MDP in Theorem \ref{thm:mdp}.
These two theorems describe the asymptotic behavior of suitably scaled
versions of the spin per site
for small values of $\alpha$.


\section{\texorpdfstring{Proof of part (a) of Theorem
\protect\ref{thm:mainsmall}}{Proof of part (a) of Theorem 4.1}}
\label{section:smallalpha}

We start by sketching how we will prove part (a) of Theorem \ref
{thm:mainsmall}.
When the quantity $\alpha$ parametrizing\vspace*{1pt} the sequence
$(\beta_n,K_n)$ satisfies $0 < \alpha< \alpha_0$,
Theorem \ref{thm:mdp} states the MDP for $S_n/n^{1 - \theta\alpha}$
under the hypotheses of Theorem \ref{thm:mainsmall} except for hypothesis
(iii$^\prime$)(b). The rate function
in this MDP is $g(x) - \inf_{y \in{\mathbb R}} g(y)$, which under
the latter hypothesis
has global minimum points at~$\pm\bar{x}$.
The MDP implies that the $P_{n,\beta_n,K_n}$-distributions of
$S_n/n^{1 - \theta
\alpha}$ put an exponentially small mass on the complement
of any open set containing the global minimum points $\pm\bar{x}$ of
the rate function.
Symmetry then yields the following weak limit, stated in Corollary \ref
{cor:weaklimit}:
\[
P_{n,\beta_n,K_n}\{S_n/n^{1 - \theta\alpha} \in dx\} \Longrightarrow
\bigl(\tfrac{1}{2}
\delta_{\bar{x}}
+ \tfrac{1}{2}\delta_{-\bar{x}} \bigr)(dx);
\]
that is, if $f$ is any bounded, continuous function, then
\[
\lim_{n \rightarrow\infty} \int_{\Lambda^n}f(S_n/n^{1-\theta
\alpha})
\,dP_{n,\beta_n,K_n}
= \int_{{\mathbb R}} f \,d \biggl(\frac{1}{2}\delta_{\bar{x}}
+ \frac{1}{2}\delta_{-\bar{x}} \biggr) = \frac{1}{2}f(\bar{x}) + \frac
{1}{2}f(-\bar{x}).
\]
In Lemma \ref{lem:uniformint} we verify that with respect to
$P_{n,\beta_n,K_n}$, the sequence $S_n/n^{1-\theta\alpha}$ is uniformly
integrable. The uniform integrability allows us
to replace the bounded, continuous function $f$ in the last display by
the absolute value function, yielding
\[
\lim_{n \rightarrow\infty} \int_{\Lambda^n} |S_n/n^{1-\theta
\alpha}|
\,dP_{n,\beta_n,K_n} =
\lim_{n \rightarrow\infty} E_{n,\beta_n,K_n} |S_n/n^{1-\theta
\alpha}|
= \bar{x}.
\]
This limit is the conclusion of part (a) of Theorem \ref{thm:mainsmall}.

We next formulate the concept of an MDP for the mean-field B--C model.
Let $(\beta_n,K_n)$ be a positive sequence converging either to a second-order
point or to the tricritical point.
Also let $\gamma$ and $u$ be real numbers satisfying $\gamma\in(0,1/2)$
and $u \in(0,1)$, and
let $\Gamma$ be a continuous
function on ${\mathbb R}$ that satisfies $\Gamma(x) \rightarrow\infty
$ as $|x|
\rightarrow\infty$.
For any subset $A$ of ${\mathbb R}$, $\Gamma(A)$ denotes the infimum of
$\Gamma$ over $A$.
We say that with respect to $P_{n,\beta_n,K_n}$, $S_n/n^{1-\gamma}$ satisfies
the MDP with
exponential speed $n^u$ and rate function $\Gamma$ if for any closed
set $F$ in ${\mathbb R}$
%
\begin{equation}
\label{eqn:uppermd}
\limsup_{n \rightarrow\infty} \frac{1}{n^{u}} \log P_{n,\beta
_n,K_n}\{
S_n/n^{1-\gamma} \in F\}
\leq-\Gamma(F)
\end{equation}
and for any open set $\Phi$ in ${\mathbb R}$
%
\begin{equation}
\label{eqn:lowermd}
\liminf_{n \rightarrow\infty} \frac{1}{n^{u}} \log P_{n,\beta
_n,K_n}\{
S_n/n^{1-\gamma} \in\Phi\}
\geq-\Gamma(\Phi).
\end{equation}
%
While an MDP is also a large deviation principle, the term MDP is often
used whenever
the exponential speed $a_n$ of the large deviation probabilities
satisfies $a_n/n \rightarrow0$
as $n \rightarrow\infty$; \cite{DemZei}, Section 3.7.

For $0 < \alpha< \alpha_0$ we now state the MDP for $S_n/n^{1-\theta
\alpha}$ with exponential speed
$n^{1 - \alpha/\alpha_0}$. The hypotheses of Theorem \ref
{thm:mainsmall} are hypotheses
(i) and (ii) of Theorem~\ref{thm:exactasymptotics}, hypotheses (iii)(a)
and (iv) of that theorem for all
$0 < \alpha< \alpha_0$, hypothesis (iii$^\prime$)(b) and the inequality $0 < \theta\alpha_0 < 1/2$.
The MDP holds under the same hypotheses except for hypothesis (iii$^\prime$)(b), which
requires that the set of global minimum points of the Ginzburg--Landau
polynomial $g$ equals $\{\pm\bar{x}\}$
for some $\bar{x}> 0$. Later in this section we will use the MDP
together with this hypothesis on the
set of global minimum points of $g$ to prove Theorem \ref{thm:mainsmall}.
Since $0 < \alpha< \alpha_0$
and $0 < \theta\alpha_0 < 1/2$, the quantities appearing in the
exponents of $n$
in the MDP satisfy $0 < \theta\alpha< 1/2$ and $0 < 1 - \alpha
/\alpha_0 < 1$.
The latter inequality implies that the exponential
speed satisfies $n^{1 - \alpha/\alpha_0} \rightarrow\infty$ as $n
\rightarrow
\infty$.
\begin{theorem}
\label{thm:mdp}
Let $(\beta_n,K_n)$ be a positive sequence parametrized by $\alpha> 0$
and converging either to a second-order point $(\beta,K(\beta))$,
$0 < \beta< \beta_c$, or to the tricritical point $(\beta_c,K(\beta
_c))$. We
assume hypotheses \textup{(i)} and
\textup{(ii)} of Theorem \ref{thm:exactasymptotics}, hypotheses
\textup{(iii)(a)} and \textup{(iv)} of that theorem
for all $0 < \alpha< \alpha_0$ and the inequality $0 < \theta\alpha
_0 < 1/2$.
Then for all $0 < \alpha< \alpha_0$, $S_n/n^{1-\theta\alpha}$
satisfies the MDP with respect to $P_{n,\beta_n,K_n}$
with exponential speed $n^{1-\alpha/\alpha_0}$ and rate function
$\Gamma(x) = g(x) - \inf_{y \in{\mathbb R}} g(y)$.
\end{theorem}

The MDP in Theorem \ref{thm:mdp} is proved exactly like the MDP in
part (a) of
Theorem~8.1 in \cite{CosEllOtt} with only changes in notation.
Rather than repeat the proof, we motivate the MDP via the related
Laplace principle.
Given $\gamma\in(0,1/2)$ and $u \in(0,1)$,
we say that with
respect to $P_{n,\beta_n,K_n}$, $S_n/n^{1-\gamma}$ satisfies the Laplace
principle with
exponential speed $n^u$ and rate function $\Gamma$ if for any bounded,
continuous function~$\psi$
\[
\lim_{n \rightarrow\infty} \frac{1}{n^{u}} \log\int_{\Lambda^n}
\exp
[n^{u}\psi(S_n/n^{1-\gamma})] \,dP_{n,\beta_n,K_n}
= \sup_{x \in{\mathbb R}}\{\psi(x) - \Gamma(x)\}.
\]
By Theorem 1.2.3 in \cite{DupEll}, if
$S_n/n^{1-\gamma}$ satisfies the Laplace principle
with exponential speed $n^u$ and rate function $\Gamma$,
then $S_n/n^{1-\gamma}$ satisfies the MDP
with the same exponential speed and the same
rate function.

Under the hypotheses of Theorem \ref{thm:mdp} we now motivate the
Laplace principle for $S_n/n^{1-\theta\alpha}$
with exponential speed $n^{1-\alpha/\alpha_0}$ and thus the MDP
stated in that theorem.
The main ideas are only sketched
because full details of the proof of an analogous Laplace principle are
given in the proof of Theorem 8.1 in \cite{CosEllOtt}.
Fix $u \in(0,1)$. If $b_n$ and $c_n$ are two positive sequences, then
we write
$b_n \asymp c_n$ if
\[
\lim_{n \rightarrow\infty} \frac{1}{n^{u}} \log b_n
= \lim_{n \rightarrow\infty} \frac{1}{n^{u}} \log c_n.
\]

We need the following lemma. It can be proved like Lemma 3.3 in
\cite{EllNew}, which applies to the Curie--Weiss
model, or like Lemma 3.2 in \cite{EllWan}, which applies to the
Curie--Weiss--Potts model.
In an equivalent form, the next lemma is well known in the literature
as the
Hubbard--Stratonovich transformation,
where it is invoked to analyze models with quadratic Hamiltonians
(see, e.g., \cite{AntRuf}, page 2363). The following lemma is also
used in the proof of Theorem \ref{thm:mainlarge}
in the next section.
\begin{lemma}
\label{lem:G}
Given a positive sequence $(\beta_n,K_n)$, let
$W_n$ be a sequence of normal
random variables with mean 0 and variance $(2\beta_n K_n)^{-1}$ defined
on a probability space $(\Omega,\mathcal{F},Q)$.
Then for any $\bar\gamma\in[0,1)$ and any bounded, continuous
function~$f$,
%
\begin{eqnarray}
\label{eqn:G}\hspace*{28pt}
&&
\int_{\Lambda^n \times\Omega}
f(S_n/n^{1-\bar\gamma} + W_n/n^{1/2-\bar\gamma})\,
d(P_{n,\beta_n,K_n} \times Q)
\nonumber\\[-8pt]\\[-8pt]
&&\qquad = \frac{1}{\int_{{\mathbb R}} \exp[-nG_{\beta
_n,K_n}(x/n^{\bar
\gamma})] \,dx} \cdot
\int_{{\mathbb R}} f(x) \exp[-nG_{\beta_n, K_n}(x/n^{\bar\gamma})]
\,dx.\nonumber
\end{eqnarray}
\end{lemma}

Let $\psi$ be any bounded, continuous function. We start our
motivation of the proof of the Laplace principle
for $S_n/n^{1-\theta\alpha}$ with exponential speed $n^{1-\alpha
/\alpha_0}$
by substituting $\bar\gamma= \theta\alpha$ and $f = \exp
(n^{1-\alpha/\alpha_0} \psi)$
into (\ref{eqn:G}), obtaining
%
\begin{eqnarray}
\label{eqn:GG}
&&
\int_{\Lambda^n \times\Omega} \exp[n^{u}\psi({S_n}/{n^{1-\gamma
}} +
{W_n}/{n^{1/2-\gamma}})] \,d(P_{n,\beta_n,K_n} \times Q)
\nonumber\\[-8pt]\\[-8pt]
&&\qquad = \frac{1}{Z_{n,\gamma}} \cdot
\int_{{\mathbb R}} \exp[n^{u}\{\psi(x) - n^{1-u}G_{\beta_n,
K_n}(x/n^{\gamma})\}] \,dx.\nonumber
\end{eqnarray}
In order to simplify the notation, we have written $\gamma$ in place
of $\theta\alpha$
and $u$ in place of $1 - \alpha/\alpha_0$. In the last display
$Z_{n,\gamma}$ is the normalization equal to
%
\begin{equation}
\label{eqn:defineznhere}
Z_{n,\gamma} = \int_{{\mathbb R}} \exp[-n G_{\beta
_n,K_n}(x/n^{\gamma})] \,dx.
\end{equation}

Let us suppose that the limit of $n^{-u}$ times the logarithm of the
right-hand side of (\ref{eqn:GG}) exists.
We then claim that since $0 < \alpha< \alpha_0$, the term
${W_n}/{n^{1/2-\gamma}}$ does not contribute to the asymptotic behavior
of the left-hand side of (\ref{eqn:GG}). From this claim it
follows that
if the limit of $n^{-u}$ times the logarithm of the right-hand side
exists, then
%
\begin{eqnarray}
\label{eqn:combinew}
&&
\int_{\Lambda^n} \exp[n^{u}\psi({S_n}/{n^{1-\gamma}}) ]
\,dP_{n,\beta_n,K_n}
\nonumber\\[-8pt]\\[-8pt]
&&\qquad \asymp\frac{1}{Z_{n,\gamma}} \cdot
\int_{{\mathbb R}} \exp[n^{u}\{\psi(x) - n^{1-u}G_{\beta_n,
K_n}(x/n^{\gamma})\}] \,dx.\nonumber
\end{eqnarray}
As on page 543 of \cite{CosEllOtt}, we justify the claim by showing
that $W_n/{n^{1/2-\gamma}}$
is superexponentially small relative to $n^{u}$ \cite{DupEll}, Theorem 1.3.3.
This holds provided $1 - 2\gamma= 1 - 2\theta\alpha> u = 1 - \alpha
/\alpha_0$,
which is valid since $0 < \theta\alpha_0 < 1/2$.
This completes our justification of the claim.

We continue our motivation of the Laplace principle for
$S_n/n^{1-\gamma}$.
The uniform convergence of $n^{1-u} G_{\beta_n,K_n}(x/n^\gamma)$ to
$g(x)$ in
hypothesis (iii)(a) of
Theorem \ref{thm:exactasymptotics} suggests that
%
\begin{eqnarray}
\label{eqn:last}
&& \int_{\Lambda^n} \exp[n^{u}\psi({S_n}/{n^{1-\gamma}}) ]
\,dP_{n,\beta_n,K_n}
\nonumber\\
&&\qquad \asymp\frac{1}{Z_{n,\gamma}} \cdot
\int_{{\mathbb R}} \exp[n^{u}\{\psi(x) - n^{1-u}G_{\beta_n,
K_n}(x/n^{\gamma})\}] \,dx \\
&&\qquad \asymp\frac{1}{\int_{{\mathbb R}} \exp[- n^{u}g(x)\}] \,dx}
\cdot
\int_{{\mathbb R}} \exp[n^{u}\{\psi(x) - g(x)\}] \,dx.\nonumber
\end{eqnarray}
The proof of this asymptotic relationship is based on hypothesis
(iii)(a) of Theorem~\ref{thm:exactasymptotics}
for $0 < \alpha< \alpha_0$,
which states that $n^{1-u} G_{\beta_n,K_n}(x/n^\gamma) = n^{\alpha
/\alpha
_0}G_{\beta_n,K_n}(x/ n^{\theta\alpha})$ converges to $g(x)$
uniformly on
compact sets,
and on several other steps, which depend in part on the lower bound in
hypothesis (iv) of Theorem
\ref{thm:exactasymptotics} for $0 < \alpha< \alpha_0$.

We define
$\bar{g} = \inf_{y \in{\mathbb R}} g(y)$. According to Laplace's method,
the asymptotic behavior of the integrals in the last line of
(\ref{eqn:last})
is governed by the maximum values of the respective integrands. Hence
%
\begin{eqnarray}
\label{eqn:why}
&&
\int_{{\mathbb R}} \exp[n^{u}\{\psi(x) - n^{1-u}G_{\beta_n,
K_n}(x/n^{\gamma})\}] \,dx \nonumber\\[-1pt]
&&\qquad \asymp\int_{{\mathbb R}} \exp[n^{u}\{\psi(x) - g(x)\}]
\,dx \\[-1pt]
&&\qquad \asymp\exp\Bigl[n^{u}\cdot\sup_{x \in{\mathbb R}}\{\psi
(x) -
g(x)\}\Bigr]\nonumber
\end{eqnarray}
and
%
\begin{eqnarray}
\label{eqn:why2}
Z_{n,\gamma} & = & \int_{{\mathbb R}} \exp[n^{u}\{- n^{1-u}G_{\beta_n,
K_n}(x/n^{\gamma})\}] \,dx \nonumber\\[-1pt]
& \asymp& \int_{{\mathbb R}} \exp[- n^{u}g(x)] \,dx \\[-1pt]
& \asymp& \exp\Bigl[- n^{u}\cdot\inf_{y \in{\mathbb R}} g(y)\Bigr]
= \exp[- n^{u}\bar{g}].\nonumber
\end{eqnarray}
Combining these two asymptotic relationships gives
\begin{eqnarray*}
&&
\int_{\Lambda^n} \exp[n^{u}\psi(S_n/n^{1-\gamma}) ] \,dP_{n,\beta
_n,K_n} \\[-1pt]
&&\qquad \asymp\frac{1}{Z_{n,\gamma}} \int_{{\mathbb R}} \exp[n^{u}\{
\psi(x) -
n^{1-u}G_{\beta_n, K_n}(x/n^{\gamma})\}] \,dx \\[-1pt]
&&\qquad \asymp\exp\Bigl[n^{u}\cdot\sup_{x \in{\mathbb R}}\bigl\{\psi(x) - \bigl(g(x) -
\bar{g}\bigr)\bigr\}\Bigr].
\end{eqnarray*}
These calculations complete the motivation that $S_n/n^{1-\gamma} =
S_n/n^{1-\theta\alpha}$
satisfies the Laplace principle and thus the MDP with exponential speed
$n^{u}= n^{1- \alpha/\alpha_0}$
and rate function $\Gamma(x) = g(x) - \bar{g}$.\vspace*{6pt}

The hypotheses of the MDP in Theorem \ref{thm:mdp} are the hypotheses
of Theorem~\ref{thm:mainsmall} except for hypothesis (iii$^\prime$)(b).
We now bring in that hypothesis, which states that the set of global
minimum points of the Ginzburg--Landau polynomial equals $\{\pm
\bar{x}\}$ for some $\bar{x}> 0$. In conjunction with the MDP we use
this hypothesis to prove part (a) of Theorem \ref {thm:mainsmall}. The
next step in that proof is contained in the following corollary, which
states that the sequence of $P_{n,\beta_n,K_n} $-distributions of
$S_n/n^{1-\theta\alpha}$ converges weakly to a symmetric sum of point
masses at $\bar{x}$ and $-\bar{x}$. This is almost immediate because up
to an additive constant the rate function in the MDP equals~$g$, and so
the $P_{n,\beta_n,K_n}$-distributions of $S_n/n^{1-\theta \alpha}$ put
an exponentially small mass on the complement of any open set
containing the global minimum points $\pm\bar{x}$ of the rate function.

We saw in the last section that the hypotheses of Theorem \ref{thm:mainsmall}
are valid for all six sequences defined in equations (\ref
{eqn:bnone})--(\ref{eqn:bnsix}) except\vspace*{1pt} for case (d)
of sequence 4, which is defined for $\ell= \ell_c$ and suitable
values of $\tilde\ell$.
As noted in the discussion leading up to Theorem \ref{thm:resultsfour},
when $\ell= \ell_c$, the set of global minimum points of $g$
equals $\{0,\pm\bar{x}\}$ for some $\bar{x}> 0$.
We are currently investigating the form of the weak limit replacing
(\ref{eqn:weaklimsmallalpha}) in the next
corollary when the set of global minimum points of $g$ has this form.
The conjecture is that in this case
there exists $0 < \lambda< 1/2$ such that
\[
P_{n,\beta_n,K_n}\{S_n/n^{1-\theta\alpha} \in dx\} \Longrightarrow
\bigl((1-2\lambda) \delta_0 + \lambda\delta_{\bar{x}} + \lambda\delta
_{-\bar{x}} \bigr) (dx).
\]
By the uniform integrability proved in Lemma \ref{lem:uniformint},
this weak limit, if true, would imply that
\[
\lim_{n \rightarrow\infty} E_{n,\beta_n,K_n}\{|S_n/n^{1 - \theta
\alpha}|\} =
2\lambda\bar{x}.
\]
\begin{cor}
\label{cor:weaklimit}
Let $(\beta_n,K_n)$ be a positive sequence parametrized by $\alpha> 0$
and converging either to a second-order point $(\beta,K(\beta))$,
$0 < \beta< \beta_c$, or to the tricritical point $(\beta_c,K(\beta
_c))$. We assume
the hypotheses of
Theorem \ref{thm:mainsmall}. Then
for all $0 < \alpha< \alpha_0$ we have the weak limit
%
\begin{equation}
\label{eqn:weaklimsmallalpha}
P_{n,\beta_n,K_n}\{S_n/n^{1-\theta\alpha} \in dx\} \Longrightarrow
\bigl(\tfrac
{1}{2} \delta_{\bar{x}} +
\tfrac{1}{2} \delta_{-\bar{x}} \bigr) (dx),
\end{equation}
where $\{\pm\bar{x}\}$ is the set of global minimum points of $g$ as specified
in hypothesis \textup{(iii$^\prime$)(b)} of Theorem \ref{thm:mainsmall}.
\end{cor}
\begin{pf} We write
$\gamma$ for $\theta\alpha$ and $u$ for $1 - \alpha/\alpha_0$.
Since $0 < \alpha< \alpha_0$, we have $0 < u < 1$, and so $n^u
\rightarrow
\infty$ as $n \rightarrow\infty$.
Let $\varepsilon> 0$ be given.
There exists $M > 0$ such that the rate function
$\Gamma(x) = g(x) - \inf_{y \in{\mathbb R}}g(y)$ in the MDP in
Theorem~\ref{thm:mdp}
is an increasing function on the interval $[M,\infty)$ and $\Gamma(M)
> 0$.
Hence the moderate deviation upper bound and symmetry imply that
\[
P_{n,\beta_n,K_n}\{|S_n/n^{1-\gamma}| \geq M \} \leq\exp[- n^u
\Gamma(M)/2].
\]
It follows that for all sufficiently large $n$, $P_{n,\beta_n,K_n}\{
S_n/n^{1-\gamma} \notin
[-2M,2M]\} \leq\varepsilon$. Thus the distributions $P_{n,\beta
_n,K_n}\{
S_n/n^{1-\gamma} \in dx\}$
are tight, and any subsequence has a weakly convergent subsubsequence
\cite{Shi}, Theorem 1, Section III.2. We now
apply the moderate deviation upper bound to any closed set $F$ in
${\mathbb R}$
not containing the global minimum points $\pm\bar{x}$ of $\Gamma$.
Since $\Gamma(F) > 0$,
we have for all sufficiently large $n$
\[
P_{n,\beta_n,K_n}\{S_n/n^{1-\gamma} \in F\} \leq\exp[- n^{u
} \Gamma
(F)/2] \rightarrow0.
\]
Thus by symmetry, any subsequence of $P_{n,\beta_n,K_n}\{
S_n/n^{1-\gamma} \in
dx\}$
has a subsubsequence converging weakly to
$(\frac{1}{2}\delta_{\bar{x}} + \frac{1}{2}\delta_{-\bar
{x}})(dx)$. This yields the weak limit in
(\ref{eqn:weaklimsmallalpha}).
\end{pf}

%

We are now ready to prove part (a) of Theorem \ref{thm:mainsmall}.
If the sequence $S_n/n^{1-\theta\alpha}$ is uniformly integrable
\cite{Bill}, Theorem 5.4, then by integrating both sides of
(\ref{eqn:weaklimsmallalpha})
with respect to the absolute value function, we obtain for all $0 <
\alpha< \alpha_0$
\[
\lim_{n \rightarrow\infty} E_{n,\beta_n, K_n} \{|S_n/n^{1-\theta
\alpha
}|\} = \bar{x}.
\]
This assertion is part (a) of Theorem \ref{thm:mainsmall}.
The required uniform integrability is proved in the next lemma from the
MDP in Theorem \ref{thm:mdp}.
\begin{lemma}
\label{lem:uniformint}
The random variables $S_n/n^{1-\theta\alpha}$ in Corollary
\ref{cor:weaklimit} are
uniformly integrable with respect to $P_{n,\beta_n,K_n}$; that is,
\[
\lim_{M \rightarrow\infty} \sup_{n \in\mathbb N}
E_{n,\beta_n,K_n} \bigl\{|S_n/n^{1-\theta\alpha}| \cdot1_{\{
|S_n/n^{1-\theta\alpha}| \geq M\}}\bigr\} = 0.
\]
\end{lemma}
\begin{pf} We write $\gamma$ for $\theta\alpha$ and $u$ for $1 -
\alpha/\alpha_0$.
$\Gamma$ denotes the rate function $g - \inf_{y \in{\mathbb R}}
g(y)$ in the
MDP in Theorem \ref{thm:mdp}.
Since $0 < \alpha< \alpha_0$, we have $0 < u < 1$, and so $n^u
\rightarrow
\infty$ as $n \rightarrow\infty$.
Let $\varepsilon> 0$ be given.
Since $g$ is a polynomial and $g(x) \rightarrow\infty$ as $|x|
\rightarrow\infty
$, $\exists M_0
\in(0,\infty)$ such that
\[
\inf_{|x| \geq M_0} \Gamma(x) \geq g(M_0)/2 > 0
\quad\mbox{and}\quad \exp\biggl[-\frac{1}{8}g(M_0)\biggr] \leq\varepsilon.
\]
The MDP in Theorem \ref{thm:mdp} implies that for all $M \geq M_0$
there exists $N_0 \in\mathbb N$
depending only on $M_0$ such that for all $n \geq N_0$
\begin{eqnarray*}
P_{n,\beta_n,K_n}\{|S_n/n^{1-\gamma}| \geq M\} & \leq&
P_{n,\beta_n,K_n}\{|S_n/n^{1-\gamma}| \geq M_0\} \\
& \leq& \exp\biggl[- \frac{1}{2}n^u \inf_{\{|x| \geq M_0\}}
\Gamma(x) \biggr].
\end{eqnarray*}
Since $|S_n| \leq n$, it follows that for all $M \geq M_0$ and for all
$n \geq N_0$
\begin{eqnarray*}
&&
E_{n,\beta_n,K_n}\bigl\{|S_n/n^{1-\gamma}| \cdot1_{\{|S_n/n^{1-\gamma}|
\geq M\}}\bigr\} \\
&&\qquad \leq
E_{n,\beta_n,K_n}\bigl\{|S_n/n^{1-\gamma}| \cdot1_{\{|S_n/n^{1-\gamma}|
\geq M_0\}}\bigr\} \\
&&\qquad \leq n^\gamma\cdot
\exp\biggl[- \frac{1}{2}n^u \inf_{|x| \geq M_0} \Gamma(x) \biggr] \\
&&\qquad \leq n^\gamma\cdot\exp\biggl[- \frac{1}{4}n^u g(M_0) \biggr].
\end{eqnarray*}
There exists $N_1 \geq N_0$ such that for all $n \geq N_1$,
$n^\gamma\cdot\exp[-\frac{1}{8}n^u g(M_0)] \leq1$.
Hence for all $M \geq M_0$ and for all $n \geq N_1$,
\[
E_{n,\beta_n,K_n}\bigl\{|S_n/n^{1-\gamma}| \cdot1_{\{|S_n/n^{1-\gamma}|
\geq M\}}\bigr\} \leq\exp\bigl[- \tfrac{1}{8}n^u g(M_0)\bigr]
\leq\exp\bigl[- \tfrac{1}{8}g(M_0)\bigr] \leq\varepsilon,
\]
which implies that for all $M \geq M_0$
\[
\sup_{n \geq N_1} E_{n,\beta_n,K_n}\bigl\{|S_n/n^{1-\gamma}| \cdot1_{\{
|S_n/n^{1-\gamma}| \geq M\}}\bigr\}
\leq\varepsilon.
\]
In addition,
\begin{eqnarray*}
&&\max_{1 \leq n < N_1} E_{n,\beta_n,K_n}\bigl\{|S_n/n^{1-\gamma}| \cdot
1_{\{|S_n/n^{1-\gamma}| \geq M\}}\bigr\} \\
&&\qquad \leq N_1^\gamma\cdot\max_{1 \leq n < N_1}
P_{n,\beta_n,K_n}\{|S_n/n^{1-\gamma}| \geq M\} \rightarrow0
\qquad\mbox{as } M
\rightarrow\infty.
\end{eqnarray*}
The last two displays complete the proof of the desired uniform
integrability. The proof of part (a) of
Theorem \ref{thm:mainsmall} is complete.
\end{pf}

In the next section we prove part (a) of Theorem \ref{thm:mainlarge}.
This theorem
gives the asymptotics of $E_{n,\beta_n,K_n}\{|S_n/n^{1-\theta\alpha
_0}|\}$ when the quantity
$\alpha$ parametrizing the sequence $(\beta_n,K_n)$ exceeds $\alpha_0$.


\section{\texorpdfstring{Proof of part (a) of Theorem
\protect\ref{thm:mainlarge}}{Proof of part (a) of Theorem 4.2}}
\label{section:largealpha}

%
Under the assumption that the quantity $\alpha$ parametrizing the
sequence $(\beta_n,K_n)$ exceeds $\alpha_0$,
part (a) of Theorem \ref{thm:mainlarge} states that
\begin{eqnarray*}
&&
\lim_{n \rightarrow\infty} n^{\theta\alpha_0} E_{n,\beta_n,K_n}\{
|S_n/n|\}
\\
&&\qquad = \lim_{n \rightarrow\infty} E_{n,\beta_n,K_n}\{
|S_n/n^{1-\theta\alpha
_0}|\} \\
&&\qquad = \bar{y}= \frac{1}{\int_{{\mathbb R}} \exp[-\tilde{g}(x)] \,dx}
\cdot
\int_{{\mathbb R}} |x| \exp[-\tilde{g}(x)] \,dx.
\end{eqnarray*}
Let $\Pi_n$ and $\Pi$ denote the probability measures on ${\mathbb R}$
defined by
%
\begin{equation}
\label{eqn:pin}\qquad
\Pi_n(dx) = \frac{1}{\int_{{\mathbb R}} \exp[-nG_{\beta_n,
K_n}(x/n^{\theta\alpha_0})] \,dx} \cdot
\exp[-nG_{\beta_n, K_n}(x/n^{\theta\alpha_0})] \,dx
\end{equation}
and
%
\begin{equation}
\label{eqn:pi}
\Pi(x) = \frac{1}{\int_{{\mathbb R}} \exp[- \tilde{g}(x)] \,dx}
\cdot
\exp[-\tilde{g}(x)] \,dx.
\end{equation}
The quantity $\bar{y}$ can be written in terms of $\Pi$ as $\int
_{{\mathbb R}}
|x| \,d\Pi$.
Part (a) of Theorem~\ref{thm:mainlarge} is proved in two steps, the
weak-convergence limit in step 1 and the uniform-integrability-type
limit in Proposition \ref{prop:weakerunifint} that yields step 2.

\textit{Step} 1. Prove that the sequence $\Pi_n$ and the sequence of
$P_{n,\beta_n,K_n}$-distributions of $S_n/n^{1-\theta\alpha_0}$
both converge weakly
to $\Pi$; that is, for any bounded, continuous function $f$,
\[
\lim_{n \rightarrow\infty} \int_{{\mathbb R}} f\, d\Pi_n = \int
_{{\mathbb R}} f\, d\Pi
\]
and
\begin{eqnarray*}
\lim_{n \rightarrow\infty} E_{n,\beta_n,K_n}\{f(S_n/n^{1-\theta
\alpha_0})\}
&=& \lim_{n \rightarrow\infty} \int_{{\mathbb R}} f(S_n/n^{1-\theta
\alpha_0})
\,dP_{n,\beta_n,K_n} \\
&=& \int_{{\mathbb R}} f\, d\Pi.
\end{eqnarray*}

\textit{Step} 2. Prove
%
\begin{eqnarray}
\label{eqn:step2}
\lim_{n \rightarrow\infty}E_{n,\beta_n,K_n}\{|S_n/n^{1-\theta
\alpha_0}|\}
& = &
\lim_{n \rightarrow\infty} \int_{{\mathbb R}} |S_n/n^{1-\theta
\alpha_0}|
\,dP_{n,\beta_n,K_n} \nonumber\\
& = & \lim_{n \rightarrow\infty} \int_{{\mathbb R}} |x|
\,d\Pi_n
= \bar{y}\\
&=& \int_{{\mathbb R}} |x|\, d\Pi.\nonumber
\end{eqnarray}
The key is to approximate the unbounded function $|x|$ by the sequence
of bounded, continuous functions
$f_j(x) = \min\{|x|,j\}$. The limits in the last display
are a consequence of the limits in
step 1 and the uniform-integrability-type limit in Proposition \ref
{prop:weakerunifint}.\vadjust{\goodbreak}

The proof of the weak-convergence limit in step 1 is given in the next theorem.
\begin{theorem}
\label{thm:scalinglargealpha}
We assume the hypotheses of Theorem \ref{thm:mainlarge}. Then
for all $\alpha> \alpha_0$ the following conclusions hold:

\textup{(a)} The sequence $\Pi_n$ defined in
(\ref{eqn:pin}) converges weakly to the probability measure $\Pi$
defined in (\ref{eqn:pi}).

\textup{(b)} The $P_{n,\beta_n,K_n}$-distributions of $S_n/n^{1-\theta
\alpha_0}$ converges weakly to $\Pi$.
\end{theorem}

The proof of this theorem relies on the following technical lemma,
which is proved in part (c) of Lemma 4.4 in \cite{CosEllOtt}.
\begin{lemma}
\label{lem:complement}
Let $(\beta_n,K_n)$ be a positive sequence parametrized by $\alpha> 0$
and converging either to a second-order point $(\beta,K(\beta))$,
$0 < \beta< \beta_c$, or to the tricritical point $(\beta,K(\beta)) =
(\beta_c,K(\beta_c))$. Assume that there exists
$\bar\gamma> 0$ and $\bar{R} > 0$ such that the sequence
\[
\xi_n = \int_{\{|x| < \bar{R}n^{\bar\gamma}\}} \exp[-nG_{\beta
_n,K_n}(x/n^{\bar\gamma})] \,dx
\]
is bounded. Then there exist constants $c_1 > 0$ and $c_2 > 0$ such
that for all sufficiently large $n$
\[
\int_{\{|x| \geq\bar{R}n^{\bar\gamma}\}} \exp[-nG_{\beta
_n,K_n}(x/n^{\bar
\gamma})] \,dx \leq c_1 \exp[-c_2 n] \rightarrow0.
\]
\end{lemma}

We now prove Theorem \ref{thm:scalinglargealpha}.

\begin{pf*}{Proof of Theorem \ref{thm:scalinglargealpha}}
We write $\gamma_0$ for $\theta\alpha_0$.
The proof follows the same pattern as the proof of Theorem 6.1 in
\cite{CosEllOtt}.
The starting point is Lemma \ref{lem:G} with $\bar{\gamma} = \gamma
_0 = \theta\alpha_0$. That lemma states
that for any bounded,
continuous function $f$
%
\begin{eqnarray}
\label{eqn:valentine}
&&
\int_{\Lambda^n \times\Omega}
f(S_n/n^{1-\gamma_0} + W_n/n^{1/2-\gamma_0})
\,d(P_{n,\beta_n,K_n} \times Q)
\nonumber\\
&&\qquad = \frac{1}{\int_{{\mathbb R}} \exp[-nG_{\beta
_n,K_n}(x/n^{\gamma_0})] \,dx} \\
&&\qquad\quad{}\times\int_{{\mathbb R}} f(x) \exp[-nG_{\beta_n, K_n}(x/n^{\gamma_0})] \,dx,\nonumber
\end{eqnarray}
where $W_n$ is a sequence of normal random variables with mean 0 and variance
$(2\beta_nK_n)^{-1}$.
Suppose that the limit of the right-hand side of (\ref
{eqn:valentine}) equals $\int_{\mathbb R} f\, d\Pi$.
Since by hypothesis $0 < \gamma_0 = \theta\alpha_0 < 1/2$,
rewriting the limit of the left-hand side of
(\ref{eqn:valentine}) in terms of characteristic
functions shows that the term $W_n/n^{1/2-\gamma_0}$ does not
contribute to this limit.
It follows that if the limit of the right-hand side of
(\ref{eqn:valentine}) equals $\int_{\mathbb R} f\, d\Pi$, then
%
\begin{eqnarray}
\label{eqn:limit2}
&&
\lim_{n \rightarrow\infty}
\int_{\Lambda^n} f(S_n/n^{1-\gamma_0}) \,dP_{n,\beta_n,K_n}\nonumber\\
&&\qquad = \lim_{n \rightarrow\infty}
\frac{1}{\int_{{\mathbb R}} \exp[-nG_{\beta_n,K_n}(x/n^{\gamma
_0})] \,dx}
\nonumber\\[-8pt]\\[-8pt]
&&\qquad\quad\hspace*{22pt}{}\times
\int_{{\mathbb R}} f(x) \exp[-nG_{\beta_n, K_n}(x/n^{\gamma_0})] \,dx
\nonumber\\
&&\qquad = \lim_{n \rightarrow\infty} \int_{{\mathbb R}} f \,
d\Pi_n=\int_{\mathbb R} f\, d\Pi.\nonumber
\end{eqnarray}

In order to calculate the limit of the sequence $\int_{{\mathbb R}} f
\,d\Pi
_n$, we appeal
to the pointwise convergence of $nG_{\beta_n, K_n}(x/n^{\gamma_0})$
to $\tilde{g}(x)$ in hypotheses (v)
of Theorem \ref{thm:mainlarge} and the lower bound in
hypotheses (iv) of Theorem \ref{thm:exactasymptotics} for $\alpha=
\alpha_0$.
This states that there exists a polynomial $H$ satisfying $H(x)
\rightarrow
\infty$ as $|x| \rightarrow\infty$
together with the following property: $\exists R > 0$ such that
$\forall n \in\mathbb N$ sufficiently large and $\forall x \in
{\mathbb R}$
satisfying $|x/n^{\gamma_0}| < R$, $n G_{\beta_n,K_n}(x/n^{\gamma
_0}) \geq H(x)$.
We then use the integrability of $\exp[-H]$
and the dominated convergence theorem to write
%
\begin{eqnarray}
\label{eqn:limtildeg}
&&\lim_{n \rightarrow\infty}
\int_{\{|x| < Rn^{\gamma_0}\}} f(x) \exp[-nG_{\beta
_n,K_n}(x/n^{\gamma_0})] \,dx\nonumber\\[-8pt]\\[-8pt]
&&\qquad =
\int_{{\mathbb R}} f(x) \exp[-\tilde{g}(x)] \,dx.\nonumber
\end{eqnarray}
In order to handle the integrals over the complementary sets $\{|x|
\geq Rn^{\gamma_0}\}$,
we appeal to Lemma \ref{lem:complement}, for which we must verify the
hypothesis.
Setting $f \equiv1$ in~(\ref{eqn:limtildeg}), we see that
\[
\lim_{n \rightarrow\infty}
\int_{\{|x| < Rn^{\gamma_0}\}} \exp[-nG_{\beta_n,K_n}(x/n^{\gamma
_0})] \,dx =
\int_{{\mathbb R}} \exp[-\tilde{g}(x)] \,dx
\]
and thus that the sequence $\int_{\{|x| < Rn^{\gamma_0}\}} \exp
[-nG_{\beta_n,K_n}(x/n^{\gamma_0})] \,dx$ is bounded.
Lemma \ref{lem:complement} with $\bar\gamma= \gamma_0$ and $\bar
{R} = R$ implies that
\[
\lim_{n \rightarrow\infty}
\int_{\{|x| \geq Rn^{\gamma_0}\}} \exp[-nG_{\beta
_n,K_n}(x/n^{\gamma_0})] \,dx = 0.
\]
It follows that
%
\begin{equation}
\label{eqn:whatever1}\quad
\lim_{n \rightarrow\infty} \int_{{\mathbb R}} f(x) \exp[-n
G_{\beta
_n,K_n}(x/n^{\gamma_0})] \,dx
= \int_{{\mathbb R}} f(x) \exp[-\tilde{g}(x)] \,dx
\end{equation}
and
%
\begin{equation}
\label{eqn:whatever2}
\lim_{n \rightarrow\infty} \int_{{\mathbb R}} \exp[-n G_{\beta
_n,K_n}(x/n^{\gamma_0})] \,dx
= \int_{{\mathbb R}} \exp[-\tilde{g}(x)] \,dx.
\end{equation}
Substituting into (\ref{eqn:limit2}) the limits in the last
two displays yields the weak convergence
asserted in parts (a) and (b) of Theorem \ref{thm:scalinglargealpha}:
\begin{eqnarray*}
\hspace*{-4pt}&&\lim_{n \rightarrow\infty} \int_{{\mathbb R}} f(S_n/n^{1-\gamma
_0}) \,dP_{n,\beta_n
,K_n} \\
\hspace*{-4pt}&&\qquad = \lim_{n \rightarrow\infty} \frac{1}{\int_{{\mathbb R}} \exp
[-n G_{\beta
_n,K_n}(x/n^{\gamma_0})] \,dx}
\cdot\int_{{\mathbb R}} f(x) \exp[-n G_{\beta_n,K_n}(x/n^{\gamma
_0})] \,dx \\
\hspace*{-4pt}&&\qquad = \lim_{n \rightarrow\infty} \int_{{\mathbb R}} f \,d\Pi_n =
\frac{1}{\int
_{{\mathbb R}} \exp[- \tilde{g}(x)] \,dx} \cdot
\int_{{\mathbb R}} f(x) \exp[- \tilde{g}(x)] \,dx = \int_{{\mathbb
R}} f \,d\Pi.
\end{eqnarray*}
The proof of the theorem is complete.
\end{pf*}

We now turn to the proof of the limit in (\ref{eqn:step2}) in
step 2, writing
$\gamma_0$ for $\theta\alpha_0$. The proof depends on an appropriate
asymptotic formula for
$E_{n,\beta_n,K_n}\{|S_n/\break n^{1-\gamma_0}|\}$, which we derive from Lemma
\ref{lem:G}. In that lemma let
$\bar\gamma= \gamma_0$, let the bounded,
continuous function $f$ equal $f_j(x) = \min\{|x|, j\}$, and send $j
\rightarrow\infty$. The monotone
convergence theorem implies that
%
\begin{eqnarray}\qquad
\label{eqn:Gextenda}
&&
\int_{\Lambda^n \times\Omega} |{S_n}/{n^{1-\gamma_0}} +
{W_n}/{n^{1/2 - \gamma_0}}|\,
d(P_{n,\beta_n,K_n} \times Q) \nonumber\\
&&\qquad = \frac{1}{\int_{{\mathbb R}} \exp[-nG_{\beta
_n,K_n}(x/n^{\gamma_0})] \,dx} \cdot
\int_{{\mathbb R}} |x| \exp[-nG_{\beta_n, K_n}(x/n^{\gamma_0})] \,dx
\\
&&\qquad = \int_{{\mathbb R}} |x| \,d\Pi_n.\nonumber
\end{eqnarray}
In this formula $W_n$
is a sequence of normal random variables with mean 0 and variance
$(2\beta_nK_n)^{-1}$ defined on a probability
space $(\Omega,\mathcal{F},Q)$, and $\Pi_n$ is the probability
measure defined in (\ref{eqn:pin}).

We write $\tilde{E}_{n,\beta_n,K_n}$ to denote expectation with
respect to
the product measure $P_{n,\beta_n,K_n}\times Q$.
Since $(\beta_n,K_n) \rightarrow(\beta,K(\beta))$, there exists a
positive constant
$c$ such that for all $n \in\mathbb N$,
\[
\tilde{E}_{n,\beta_n,K_n}\{|W_n/n^{1/2 - \gamma_0}|\} \leq c/n^{1/2
- \gamma_0}.
\]
Thus
\begin{eqnarray*}
&&
\tilde{E}_{n,\beta_n,K_n}\{|S_n/n^{1-\gamma_0} + W_n/n^{1/2 - \gamma
_0}| \} +
c/n^{1/2 - \gamma_0} \\
&&\qquad \geq E_n\{|S_n/n^{1 - \gamma_0}|\} \geq
\tilde{E}_{n,\beta_n,K_n}\{ |S_n/n^{1 - \gamma_0} + W_n/n^{1/2 -
\gamma_0}| \} -
c/n^{1/2 - \gamma_0}.
\end{eqnarray*}
Suppose that we could prove
\[
\lim_{n \rightarrow\infty} \int_{{\mathbb R}} |x| \,d\Pi_n = \bar
{y}= \int_{{\mathbb R}}
|x| \,d\Pi.
\]
Since $0 < \gamma_0 = \theta\alpha_0 < 1/2$,
we would then obtain from (\ref{eqn:Gextenda}) the desired limit
%
\begin{eqnarray}
\label{eqn:limitrefineda}
&&
\lim_{n \rightarrow\infty} E_{n,\beta_n,K_n}\{|S_n/n^{1 - \gamma
_0}| \}\nonumber \\
&&\qquad = \lim_{n \rightarrow\infty}
\tilde{E}_{n,\beta_n,K_n}\{|S_n/n^{1 - \gamma_0} + W_n/n^{1/2 -
\gamma_0}| \} \\
& &\qquad = \lim_{n \rightarrow\infty}
\int_{{\mathbb R}} |x| \,d\Pi_n = \bar{y}= \int_{{\mathbb R}} |x|
\,d\Pi.\nonumber
\end{eqnarray}

We complete the proof of part (a) of Theorem \ref{thm:mainlarge} by
showing the limit in the last line of (\ref{eqn:limitrefineda}).
Part (a) of Theorem \ref{thm:scalinglargealpha} shows that the
sequence $\Pi_n$ converges weakly to $\Pi$.
According to a standard result, the limit in the last line of
(\ref{eqn:limitrefineda})
would follow immediately from the weak convergence of $\Pi_n$ to $\Pi$
if we could prove that $\Pi_n$ satisfies the following
uniform-integrability estimate:
\[
\lim_{j \rightarrow\infty} \sup_{n \in\mathbb N} \int_{\{|x| > j\}
} |x| \,d\Pi
_n = 0.
\]
The next proposition shows that the limit in the last line of
(\ref{eqn:limitrefineda})
is a consequence of a condition that is weaker than uniform integrability.
\begin{prop}
\label{prop:weakerunifint}
Let $\tilde{\Pi}_n$ be a sequence of probability measures on
${\mathbb R}$ that
converges weakly to a probability measure
$\tilde{\Pi}$ on ${\mathbb R}$. Assume in addition that $\int
_{{\mathbb R}} |x| \,d\tilde{\Pi}<
\infty$ and that
%
\begin{equation}
\label{eqn:weakerunifint}
\lim_{j \rightarrow\infty} \limsup_{n \rightarrow\infty} \int_{\{
|x| > j\}}
|x| \,d\tilde{\Pi}_n = 0.
\end{equation}
It then follows that
\[
\lim_{n \rightarrow\infty} \int_{{\mathbb R}} |x| \,d\tilde{\Pi}_n
= \int_{{\mathbb R}} |x|
\,d\tilde{\Pi}.
\]
\end{prop}
\begin{pf} For $j \in\mathbb N$, $f_j$ denotes the bounded, continuous
function that equals $|x|$ for
$|x| \leq j$ and equals $j$ for $|x| > j$. Then
\begin{eqnarray*}
&&
\biggl|\int_{{\mathbb R}} |x| \,d\tilde{\Pi}_n - \int_{{\mathbb R}} |x|
\,d\tilde{\Pi}\biggr| \\
&&\qquad \leq
\int_{{\mathbb R}} \bigl| |x| - f_j \bigr| \,d\tilde{\Pi}_n +
\biggl|\int_{{\mathbb R}} f_j \,d\tilde{\Pi}_n - \int_{{\mathbb R}} f_j
\,d\tilde{\Pi}\biggr| + \int_{{\mathbb R}} \bigl| |x|
- f_j \bigr|\, d\tilde{\Pi}\\
&&\qquad \leq
2 \int_{\{|x| > j\}} |x| \,d\tilde{\Pi}_n +
\biggl|\int_{{\mathbb R}} f_j \,d\tilde{\Pi}_n - \int_{{\mathbb R}} f_j
\,d\tilde{\Pi}\biggr| + 2 \int_{\{|x| >
j\}} |x| \,d\tilde{\Pi}.
\end{eqnarray*}
Since $\tilde{\Pi}_n \Rightarrow\tilde{\Pi}$, we have $\int
_{{\mathbb R}} f_j \,d\tilde{\Pi}_n
\rightarrow\int_{{\mathbb R}} f_j \,d\tilde{\Pi}$, and therefore
\begin{eqnarray*}
&&
\limsup_{n \rightarrow\infty} \biggl|\int_{{\mathbb R}} |x| \,d\tilde{\Pi
}_n - \int_{{\mathbb R}} |x|
\,d\tilde{\Pi}\biggr| \\
&&\qquad \leq
2 \limsup_{n \rightarrow\infty} \int_{\{|x| > j\}} |x| \,d\tilde{\Pi
}_n + 2 \int
_{\{|x| > j\}} |x| \,d\tilde{\Pi}.
\end{eqnarray*}
By the assumptions on $\tilde{\Pi}_n$ and $\tilde{\Pi}$, both terms
on the
right-hand side of this inequality
converge to 0 as $j \rightarrow\infty$. This completes the proof.
\end{pf}

In order to justify the limit in the last line of
(\ref{eqn:limitrefineda}),
we must verify the hypotheses of Proposition \ref{prop:weakerunifint}
for the measures $\Pi_n$ and $\Pi$
defined in (\ref{eqn:pin}) and (\ref{eqn:pi}).
Clearly the measure $\Pi$ defined in (\ref{eqn:pi})
satisfies $\int_{{\mathbb R}} |x| \,d\Pi< \infty$.
We now verify the condition in (\ref{eqn:weakerunifint}) for
the measures $\Pi_n$
defined in (\ref{eqn:pin}). For any $j \in\mathbb N$ and all
sufficiently large $n$
we will find quantities $A_j$, $B_n$ and $C_n$ with the following properties:
\[
\int_{\{|x| > j\}} |x| \,d\Pi_n \leq A_j + B_n + C_n,
\]
$A_j \rightarrow0$ as $j \rightarrow\infty$,
$B_n \rightarrow0$ as $n \rightarrow\infty$ and $C_n \rightarrow0$
as $n \rightarrow\infty$.
It follows from these properties that
%
\begin{equation}
\label{eqn:finisha}\qquad
\lim_{j \rightarrow\infty} \limsup_{n \rightarrow\infty}
\int_{\{|x| > j\}} |x| \,d\Pi_n \leq\lim_{j \rightarrow\infty} A_j +
\lim_{n \rightarrow\infty} B_n + \lim_{n \rightarrow\infty} C_n = 0.
\end{equation}
This yields the limit in (\ref{eqn:weakerunifint}), proving
step 2 and thus completing the proof
of part~(a) of Theorem \ref{thm:mainlarge}.

We now specify the quantities $A_j$, $B_n$, and $C_n$ having the
properties in the preceding paragraph.
Given positive integers $j$ and $n$, let $R$ and $c$ be positive
numbers that satisfy
$c > R$ and that will be specified below. We then partition
the set $\{|x| > j\}$ into the following three subsets:
\begin{eqnarray*}
\{|x| > j\} & = & [\{|x| > j\} \cap\{|x/n^{\gamma_0}| < R\} ] \\
& &{} \cup[\{|x| > j\} \cap\{R \leq|x/n^{\gamma_0}| < c\} ] \cup
[\{|x| > j\} \cap\{|x/n^{\gamma_0}| \geq c\} ].
\end{eqnarray*}
Since for all $n$
\[
\{|x| > j\} \subset[\{|x| > j\} \cap\{|x/n^{\gamma_0}| < R\} ] \cup
\{R \leq|x/n^{\gamma_0}| < c\} \cup\{|x/n^{\gamma_0}| \geq c\},
\]
it follows that for all $n$
%
\begin{eqnarray}
\label{eqn:threeintegrals}\qquad\quad
\int_{\{|x| > j\}} |x| \,d\Pi_n & \leq& \int_{\{|x| > j\} \cap\{
|x/n^{\gamma_0}| < R\}} |x| \,d\Pi_n +
\int_{\{R \leq|x/n^{\gamma_0}| < c\}} |x| \,d\Pi_n\nonumber\\[-8pt]\\[-8pt]
&&  + \int_{\{
|x/n^{\gamma_0}| \geq c\}} |x| \,d\Pi_n.\nonumber
\end{eqnarray}

We next estimate each of these three integrals.
The convergence proved in (\ref{eqn:whatever2}) implies that
the sequence
$1/Z_n$ is positive and bounded. By hypothesis (iv) of Theorem \ref
{thm:exactasymptotics} for $\alpha= \alpha_0$
there exists $R > 0$ such that for all sufficiently large $n \in
\mathbb N$
and all $x \in{\mathbb R}$ satisfying $|x/n^{\gamma_0}| < R$
\[
n G_{\beta_n,K_n}(x/n^{\gamma_0}) \geq H(x),
\]
where $H$ is a polynomial satisfying $H(x) \rightarrow\infty$ as $|x|
\rightarrow
\infty$.
Since $\exp[- H(x)]$ is integrable, for all sufficiently large $n$ we
estimate the first integral on the right-hand side
of (\ref{eqn:threeintegrals}) by
%
\begin{eqnarray}
\label{eqn:Aj}\qquad
&&
\int_{\{|x| > j\} \cap\{|x/n^{\gamma_0}| < R\}} |x| \,d\Pi_n
\nonumber\\[-8pt]\\[-8pt]
&&\qquad \leq A_j =
\mathrm{const} \cdot\int_{\{|x| > j\}} |x| \exp[-H(x)] \,dx \rightarrow0
\qquad\mbox{as } j \rightarrow\infty.\nonumber
\end{eqnarray}

By part (a) of Lemma 4.4 in \cite{CosEllOtt}, there exists $c > 0$ and
$D > 0$
such that $G_{\beta_n,K_n}(x)
\geq D x^2$ for all $|x| \geq c$; thus for all $n \in\mathbb N$
and all $x \in{\mathbb R}$ satisfying $|x/n^{\gamma_0}| \geq c$,
\[
n G_{\beta_n,K_n}(x/n^{\gamma_0}) \geq n D (x/n^{\gamma_0})^2 = n^{1 -
2\gamma_0} D x^2.
\]
Without loss of generality $c$ can be chosen to be larger
than the quantity $R$ specified in the preceding paragraph.
Since the sequence $1/Z_n$ is bounded,
we estimate the third integral on the right-hand side of
(\ref{eqn:threeintegrals}) by
%
\begin{eqnarray}
\label{eqn:Cn}
&&
\int_{\{|x/n^{\gamma_0}| \geq c\}} |x| \,d\Pi_n \nonumber\\
&&\qquad \leq C_n = \frac{1}{Z_n} \cdot\int_{\{|x/n^{\gamma_0}|
\geq c\}} |x| \exp[- n^{1-2\gamma_0} D x^2] \,dx \\
&&\qquad \leq
\mathrm{const} \cdot n^{2\gamma_0 - 1} \cdot\exp[- n c^2 D]
\rightarrow0\qquad
\mbox{as } n \rightarrow\infty.\nonumber
\end{eqnarray}

With these choices of $R$ and $c$ we estimate the middle integral on
the right-hand side of
(\ref{eqn:threeintegrals}) by
\begin{eqnarray*}
&&
\int_{\{R \leq|x/n^{\gamma_0}| < c\}} |x| \,d\Pi_n \\
&&\qquad \leq B_n = c n^{\gamma_0} \cdot\Pi_n\{|x/n^{\gamma_0}| \geq R\}
\\
&&\qquad = c n^{\gamma_0} \cdot\frac{1}{Z_n} \cdot\int_{\{|x/n^{\gamma
_0}| \geq R\}} \exp[-n G_{\beta_n,K_n}(x/n^{\gamma_0})] \,dx.
\end{eqnarray*}
Since the sequence $1/Z_n$ is bounded, the display after (\ref{eqn:limtildeg})
and Lemma \ref{lem:complement} with
$\bar{R} = R$ and $\bar\gamma= \gamma_0$ imply
the existence of constants $c_1 > 0$ and $c_2 > 0$ such that
%
\begin{eqnarray}
\label{eqn:bn}
\int_{\{R \leq|x/n^{\gamma_0}| < c\}} |x| \,d\Pi_n \leq B_n \leq c
n^{\gamma_0} \cdot\mathrm{const } \cdot c_1 \exp[-c_2 n]
\rightarrow0 \nonumber\\[-8pt]\\[-8pt]
\eqntext{\mbox{as } n \rightarrow\infty.}
\end{eqnarray}

Together, equations (\ref{eqn:Aj}), (\ref{eqn:bn}) and (\ref
{eqn:Cn}) prove (\ref{eqn:finisha}), which completes the proof of part (a)
of Theorem \ref{thm:mainlarge}.\qed


%
\printaddresses

\end{document}